\documentclass{amsart}
\usepackage{xypic}

\newtheorem{thm}{Theorem}[section]
\newtheorem*{thm*}{Theorem}

\newtheorem{cor}[thm]{Corollary}
\newtheorem{dfn}[thm]{Definition}
\newtheorem{ex}[thm]{Example}
\newtheorem{lemma}[thm]{Lemma}

\newtheorem{prop}[thm]{Proposition}

\newtheorem{remark}[thm]{Remark}

\renewcommand{\cH}{{\mathcal H}}
\newcommand{\cG}{{\mathcal G}}

\renewcommand{\cL}{{\mathcal L}}

\newcommand{\cU}{{\mathcal U}}

\newcommand{\cW}{{\mathcal W}}

\newcommand{\E}{{\mathcal E}}

\newcommand{\F}{{\mathcal F}}

\newcommand{\mR}{{\mathbb R}}

\newcommand{\oF}{\overline{{\mathcal F}}}

\newcommand{\oL}{\overline{L}}
\newcommand{\Oq}{{\mathbf O(q)}}

\newcommand{\SOq}{{\mathbf{SO}(q)}}
\newcommand{\SO}{{\mathbf{SO}}}
\newcommand{\soq}{{{\mathfrak{so}}(q)} }
\newcommand{\bas}{{\rm bas}\, }

\newcommand{\codim}{{\operatorname{codim\, }}}

\newcommand{\whF}{{\widehat{\F}}}

\newcommand{\whL}{{\widehat L}}
\newcommand{\whM}{{\widehat M}}

\newcommand{\whW}{{\widehat W}}

\newcommand{\mfk}{{\mathfrak k}}
\newcommand{\mfg}{{\mathfrak g}}

\newcommand{\umfg}{\underline{\mathfrak g}}
\newcommand{\mfa}{{\mathfrak a}}

\newcommand{\Iso}{{\operatorname{Iso}}}

\newcommand{\Diff}{{\operatorname{Diff}}}

\newcommand{\ZZ}{{\mathbb  Z}}
\newcommand{\NN}{{\mathbb  N}}
\newcommand{\RR}{{\mathbb  R}}

\begin{document}

\title{Localization of Basic Characteristic Classes}

\author{Dirk T\"oben}
 \thanks{The author was supported  by the Schwerpunktprogramm SPP 1154  of the DFG}
\address{Dirk T\"oben, Mathematisches Institut, Universit\"at zu K\"oln, Weyerthal 86-90, 50931 K\"oln, Germany}
\email{dtoeben@math.uni-koeln.de}




\begin{abstract}
We introduce basic characteristic classes and numbers as new invariants for Riemannian foliations. If the ambient Riemannian manifold $M$ is complete, simply connected (or more generally if the foliation is a transversely orientable Killing foliation) and if the space of leaf closures is compact, then the basic characteristic numbers are determined by the infinitesimal dynamical behavior of the foliation at the union of its closed leaves. In fact, they can be computed with an Atiyah-Bott-Berline-Vergne-type localization theorem for equivariant basic cohomology.
\end{abstract}

\maketitle

\tableofcontents


\section{Introduction}
The Poincar\'e-Hopf Theorem states that by properly counting the singularities of a vector field, or more precisely by adding their indices, one obtains the Euler characteristic. In modern terms this result can be viewed as the localization of the Euler class, where the indices are residual data; a proof of the Poincar\'e-Hopf Theorem with this method was first given by Chern for dimension two. For a Killing field, an infinitesimal isometric motion, on a Riemannian manifold Bott (\cite{Bott}) localized polynomials of top degree in the Pontryagin classes of the manifold to its singularities. Later Berline and Vergne on the one hand (\cite{BerlineVergne}) and Atiyah and Bott on the other hand reproved this result using equivariant cohomology. Atiyah and Bott used their localization formula for equivariant cohomology of torus actions (\cite[Sections 3 and 8]{AtiyahBott}). We aim to define basic characteristic numbers of Riemannian foliations and compute them in the same spirit by using equivariant basic cohomology introduced in \cite{GT2010}.

Let $(M,\F)$ be a transversely oriented Riemannian foliation of codimension $q$, i.e.~a foliation with locally equidistant leaves (for the precise definition see Section \ref{sec:Riemfol}). The characteristic forms of the normal bundle $\nu\F$ of $\F$ with respect to a Riemannian basic connection are basic forms. The corresponding cohomology classes in the basic cohomology ring $H^*(M,F)$ are called basic characteristic classes, and they form a subring which we call the basic Euler-Pontryagin ring (Section \ref{sec:BasicCharacteristicClasses}). Via a transverse version of integration $\int_\F:H^q(M,\F)\stackrel{\cong}{\to}\RR$ constructed by Sergiescu (Section \ref{sec:PoincareDuality}), we define the basic characteristic number $\int_\F p$ for an element $p$ of the basic Euler-Pontryagin ring of degree $q$. The set of these numbers is a new invariant for Riemannian foliations. We emphasize a difference of our basic point of view to other papers as e.g.~\cite{lazarovpasternack}. There the characteristic forms  of the normal bundle of the foliation are considered in $H^*(M)$ instead of $H^*(M,\F)$ which leads under the additional assumption that the normal bundle is trivial to secondary characteristic classes.

The Molino sheaf of a Riemannian foliation gives rise (in the simply-connected case or more generally for Killing foliations) to a transverse action of an abelian Lie algebra $\mfa$, called structural Killing algebra, whose orbits of the leaves are the leaf closures (Section \ref{sec:transverseaction}). In this respect this action describes the transverse dynamical behavior of the foliation. According to \cite{GT2010} a transverse action allows one to construct the equivariant basic cohomology algebra of $(M,\F)$ (see Section \ref{sec:EquivariantBasicCohomology}) which we expect to link cohomological and dynamical data. We extend the basic characteristic class $p(\Omega)$, where $\Omega$ is the curvature $2$-form of the canonical Riemannian basic connection (which is $\mfa$-invariant) to an equivariantly basic characteristic class $p(\Omega+L_X)$ in analogy to \cite{BerlineVergne} and \cite[Section 8]{AtiyahBott} and localize it via an Atiyah-Bott-Berline-Vergne type localization formula for equivariant basic cohomology to the union $C$ of closed leaves (see below). Our main result is that the characteristic number $\int_\F p$ is determined by the infinitesimal behavior of the foliation at $C$.
\begin{thm*}[ABBV-type localization formula]
For any closed equivariantly basic form $\omega$ we have
$$
\int_\F\omega=\sum_i c_i\int_{C_i/\F}\left(\frac{i^*_{C_i}\omega}{e_\mfa(\nu C_i,\whF)}\right).
$$
\end{thm*}
Here the $C_i$ are the connected components of $C$ and $e_\mfa(\nu C_i,\whF)$ is the equivariant basic Euler class of the foliated normal bundle $\nu C$ whose foliation is induced from the natural foliation $\whF$ on $\nu\F$. Note that a basic form of $(C_i,\F)$ descends to a form on the orbifold $C_i/\F$. The integral $\int_{C_i/\F}$ denotes standard integration on the orbifold $C_i/\F$, and the $c_i$ are constants (see Section \ref{sec:PoincareDuality}).

\begin{thm*}
Let $\F$ be a Riemannian foliation of even codimension $q$ on a connected, complete, simply-connected Riemannian manifold $M$ such that $M/\oF$ is compact. Moreover, assume that the closed leaves $\{L_i\}$ are isolated. Then
$$
\int_\F p(\Omega)=\sum_i (-1)^{q/2} c_i\frac{p(L_X)}{\prod_{\alpha\in\Delta_i} \alpha(X)}.
$$
for a basic characteristic class $p$ of degree $q$.
Here $\Delta_i$ is the set of (nontrivial) roots of the isotropy representation of the transverse action at $\nu_x C_i$ for any point $x\in C_i$.
\end{thm*}
The corollary is more generally true for transversely oriented Killing foliations on complete Riemannian manifolds with $M/\oF$ compact.

A crucial technical tool that is interesting in its own right is a basic version of the Thom isomorphism for foliated bundles over Riemannian foliations. We prove it in greater generality than needed here for possible further applications.

In Section \ref{sec:foliatedbundles} we define foliated bundles and prove invariance of basic cohomology under basic homotopies.

Then we briefly remind the reader of Molino's structure theory of Riemannian foliations and recall how the Molino sheaf gives rise to a transverse action for certain Riemannian foliations, namely Killing foliations, in Section \ref{sec:Molino}.

In Section \ref{sec:PoincareDuality} we recall transverse integration as defined by Sergiescu which he used to prove basic Poincar\'e duality for Riemannian foliations. We will need this integration later to define basic characteristic numbers.

In Section \ref{sec:BasicThom} we prove a basic version of the Thom isomorphism for foliated bundles over Riemannian foliations. The main application in this paper is the naturally foliated normal bundle $(\nu C,\whF)$, where $C$ is the union of closed leaves of a Riemannian foliation. For a Killing foliation we extend the basic Thom isomorphism to an equivariant basic Thom isomorphism in Section \ref{sec:EquivariantBasicCohomology} after recalling the notion of equivariant basic cohomology.

In Section \ref{sec:localization} we define basic characteristic classes for the normal bundle of a Riemannian foliation and we introduce basic characteristic numbers as new invariants for Riemannian foliations. We extend the basic characteristic classes to equivariant basic characteristic classes for Killing foliations. We then prove an Atiyah-Bott-Berline-Vergne type localization formula in equivariant basic cohomology (see above) using the equivariant basic Thom isomorphism for $(\nu C,\whF)$. The application of this localization formula to equivariant basic characteric classes finally allows us to compute the corresponding basic characteristic numbers in local terms around $C$, see the above theorem.

In the appendix we prove that a Riemannian foliation on a manifold $M$, which is complete with respect to a compatible bundle-like metric, and for which $M/\oF$ compact (e.g.~if $M$ is compact) has a good saturated cover. The latter notion is needed for the proof of the basic Thom isomorphism for Riemannian foliations.

{\em Acknowledgements:} The author is grateful to Oliver Goertsches and to Steven Hurder for various very helpful discussions.

\section{Foliated bundles}\label{sec:foliatedbundles}
Let $(M,\F)$, $(X,\cG)$ and $(E,\E)$ be foliated manifolds and let $\pi:E\to M$ be a surjective map. In this paper all maps are assumed to be smooth. Then $\pi$ is called a {\em foliated bundle} with typical fiber $X$, if the following holds: there is a cover $\{U_\alpha\}$ of $M$ and diffeomorphisms
$$
\phi:(E|{U_\alpha},\E)=(\pi^{-1}(U_\alpha),\E)\to (U_\alpha,\F)\times (X,\{*\}),
$$
of foliations, the trivializations, such that
$\pi\circ \phi^{-1}:U_\alpha\times X\to U_\alpha$ is the projection. That means, with respect to trivializations, $(E,\E)$ is locally the product foliation of $(U_\alpha,\F)$ and the trivial foliation of $X$ by points.

For $\alpha,\beta$ the map
$$
\phi_\alpha\circ\phi_\beta^{-1}:(U_{\alpha}\cap U_{\beta})\times X \to (U_{\alpha}\cap U_{\beta})\times X
$$
restricted to each fiber is an element of the group of diffeomorphisms $\Diff(X)$ of $X$. We obtain {\em transition maps}
$$
\gamma_{\alpha\beta}:U_\alpha\cap U_\beta\to \Diff(X)
$$
$$
\gamma_{\alpha\beta}(p)=\phi_\alpha\circ\phi_\beta^{-1}|\{p\}\times X
$$
that are constant along the leaves of $\F|(U_\alpha\cap U_\beta)$ and which fulfill the cocycle condition $\gamma_{\alpha\beta}(p)\circ\gamma_{\beta\delta}(p)=\gamma_{\alpha\delta}(p)$ for all $p\in U_\alpha\cap U_\beta \cap U_\delta$.

Conversely, a manifold $X$, a foliated manifold $(M,\F)$, a covering $\{U_\alpha\}$ of $M$, and transition maps $\{\gamma_{\alpha\beta}\}$ that are constant along the leaves of $\F|(U_\alpha\cap U_\beta)$ yield a foliated bundle $\pi:(E,\E)\to (M,\F)$ with typical fiber $X$. If moreover the values of the transition maps lie in a subgroup $G$ of $\Diff(X)$ we say that the structure group of the foliated bundle is $G$. If $X$ is a vector space and $G=\mathbf{GL}(X)$ we call $(E,\E)$ a {\em foliated vector bundle}.

A map $f:(M,\F)\to (M',\F')$ between foliated manifolds is called {\em foliate}, if it respects foliation, i.e., maps leaves to leaves. A homotopy $H:(M,\F)\times [0,1]\to (M',\F')$ is {\em basic} if $H_t:=H(\,\cdot\,,t):(M,\F)\to (M',\F')$ is foliate for all $t\in [0,1]$. For a foliated manifold $(M,\F)$ we consider the complex
$$
\Omega^*(M,\F)=\{\omega\in\Omega^*(M)\mid \iota_X\omega=0, \cL_X\omega=0\ \mbox{for all}\ X\in\Xi(\F) \}
$$
of basic forms; here $\Xi(\F)$ is the space of vector fields tangential to $\F$. It is a differential subcomplex of $\Omega^*(M)$  whose cohomology we denote by $H^*(M,\F)$. This is the basic cohomology ring of $(M,\F)$. The pull back $f^*$ of a foliate map $f:(M,\F)\to (M',\F')$ respects the basic complexes and therefore induces a homomorphism $f^*:H^*(M',\F')\to H^*(M,\F)$.

\begin{prop}\label{prop:simpleKunneth}
Let $(M,\F)$ be a foliated manifold and $\pi:(M\times\RR^n,\F\times\{*\})\to (M,\F)$ the trivial foliated bundle over it.
Then $H^*(M\times \RR^n,\F\times\{*\})$ and $H^*(M,\F)$ are isomorphic via $\pi^*$ which has the inverse map $s^*$, where $s:M\to M\times \RR^n$ is the zero section.
\end{prop}
\begin{proof}
By induction on $n$ it is sufficient to consider the foliated bundle $\pi:(M\times\RR,\F\times\{*\})\to (M,\F)$. For the moment we forget the foliated structures and regard $\pi:M\times \RR\to M$ as an ordinary bundle.
For a chart domain $U\subset M$ every ordinary form on $U\times\RR$ is uniquely a linear combination of the two types of forms $(\pi^*\phi)f(x,t)$ and $(\pi^*\phi)f(x,t)dt$, where $\phi$ is a form on the base $U$ and $f$ a function on $U\times\RR$.
The operator $K:\Omega^*(M\times\RR)\to\Omega^{*-1}(M\times\RR)$ is defined locally by
\begin{align*}
&(\pi^*\phi)f(x,t)\mapsto 0,\\
&(\pi^*\phi)f(x,t)dt\mapsto (\pi^*\phi)\int_0^tf.
\end{align*}
It is well-defined. Now we take the foliated structures into consideration. With respect to trivializations of $\pi$ over foliated chart domains of $(M,\F)$ the restriction of the map $\pi$ looks like the projection $\rho:\RR^p\times\RR^q\times\RR\to\RR^p\times\RR^q$ where the foliations have leaves of the form $\RR^p\times \{*\}$. Now every every basic form on $\RR^{p+q+1}$ is uniquely a linear combination of the two types of forms $(\rho^*\phi)f(x,y,t)$ and $(\rho^*\phi)f(x,y,t)dt$, with $\phi$ an ordinary form on the base $\RR^{p+q}$ and $f$ a function on $\RR^{p+q+1}$ which is independent of $x$. We see that $K$ maps basic forms to basic forms and we denote its restriction by $K_b:\Omega^*(M\times\RR,\F\times\{*\})\to\Omega^{*-1}(M\times\RR,\F\times\{*\})$.
It is known (see e.g.~\cite[p.~34]{BottTu}) that $K$ is a homotopy operator with $1-\pi^*s^*=(-1)^k(dK-Kd)$ on $\Omega^k(M\times\RR)$.  Since $\pi$ and $s$ are foliate the same identity holds for $K_b$ on the basic subcomplex. Thus the map $s^*:H^*(M\times\RR,\F\times\{*\})\to H^*(M,\F)$ and its inverse map $\pi^*$ are isomorphisms.
\end{proof}
\begin{cor}\label{cor:homotopicmaps}
Let $f_0,f_1:(M,\F)\to (M',\F')$ be foliate maps that are basically homotopic via a foliate map $f:M\times [0,1]\to M'$. Then $f_0^*=f_1^*:H^*(M',\F')\to H^*(M,\F)$.
\end{cor}
\begin{proof}
As in \cite[Corollary 4.1.2]{BottTu}. For $i=0,1$ let $s_i:M\to M\times \{i\}\subset M\times [0,1]$ be the natural inclusion. Since both $s_0^*$ and $s_1^*$ invert $\pi^*$ in basic cohomology, by Proposition \ref{prop:simpleKunneth} they are identical. Thus $f_0^*=(F\circ s_0^*)=s_0^*\circ F^*=s_1^*\circ F^*=(F\circ s_1^*)=f^*_1$.
\end{proof}
This implies that two foliated manifolds that are equivalent up to (differentiable) basic homotopies have isomorphic basic cohomology rings. In particular this is true for diffeomorphism of foliated manifolds. It is interesting that basic cohomology is not an invariant of homeomorphisms, as shown by El Kacimi and Nicolau, see \cite[p.~628]{AlNi}. In fact their main result of is that a homeomorphism $h:(M,\F)\to (M',\F')$ between two {\em Riemannian} foliations still yields an isomorphism $h^*:H^*(M',\F')\to H^*(M,\F)$ of basic cohomology.

We will also need the homotopy invariance of the compact vertical basic cohomology of foliated bundles. For a foliated bundle  $\pi:(E,\E)\to (M,\F)$ let $\Omega_{cv}(E,\E)$ be the complex of basic forms on $E$ with compact vertical support.
\begin{prop}\label{prop:cohomologybundles}
Let $(E,\E)$ be a foliated bundle over some foliated manifold $(M,\F)$ and $(E\times\RR,\E\times\{*\})$ be the corresponding foliated product bundle over $(M\times \RR,\F\times\{*\})$. Let $\pi:E\times\RR\to E$ be the projection.
Then $H^*_{cv}(E\times \RR,\E\times\{*\})$ and $H^*_{cv}(E,\E)$ are isomorphic via $\pi^*$ which has the inverse map $s^*$, where $s:E\to E\times \RR; v\mapsto (v,0)$.
\end{prop}
\begin{proof}
Basically we replace $(M,\F)$ in the proof of Proposition \ref{prop:simpleKunneth} with $(E,\E)$. We consider a trivialization of $(E,\E)$ over a foliation chart domain $U\subset M$. Every basic form on $E|U\times\RR$ is uniquely a linear combination of the two types of forms $(\pi^*\phi)f(v,t)$ and $(\pi^*\phi)f(v,t)dt$, where $\phi\in\Omega^*(E|U,\E)$ and $f\in \Omega^0_{cv}(E|U\times\RR,\E\times \{*\})$, i.e.~, a basic function with compact vertical support. Note that the operator $K$ defined as in the mentioned proof is a map $K:\Omega^*_{cv}(E\times\RR,\E\times\{*\})\to\Omega^{*-1}_{cv}(E\times\RR,\E\times\{*\})$. Following the computation in \cite[p.~34]{BottTu}) we deduce that $K$ is a homotopy operator with $1-\pi^*s^*=(-1)^k(dK-Kd)$ on $\Omega^k(M\times\RR)$.
\end{proof}
\begin{cor}\label{cor:homotopicbundles}
Let $\pi:(E,\E)\to (M,\F)$ and $\pi':(E',\E')\to (M',\F')$ be foliated bundles. Moreover, let $f:(M\times \RR,\F\times \{*\})\to (M',\F')$ and $\overline f:(E\times \RR,\E\times \{*\})\to (E',\E')$ be foliate maps such that
$$
\xymatrix{
  (E\times \RR,\E\times \{*\})  \ar[r]^{\overline f}\ar[d]_{\pi\times {\rm id}_\RR}& (E',\E')\ar[d]^{\pi'} \\
  (M\times \RR,\F\times \{*\}) \ar[r]^f& (M',\F').}
$$
is commutative. Then $\overline f_0^*=\overline f_1^*:H_{cv}^*(E',\E')\to H_{cv}^*(E,\E)$.
\end{cor}
\begin{proof}
As in the proof of Corollary \ref{cor:homotopicmaps}.
\end{proof}
We will later (in the proof of Theorem \ref{thm:thomiso}) use this result in the following situation. Let $f_0,f_1:(M,\F)\to (M',\F')$ be foliate maps that are basically homotopic via a foliate map $f:(M\times [0,1],\F\times \{*\})\to (M',\F')$. Moreover, let $\pi':(E',\E')\to (M',\F')$ be a foliated bundle with a basic connection, that means a connection that is basic (horizontal and invariant) with respect to the horizontal lifts of all $\F'$-tangential vector fields. Its parallel translation maps respect the foliation $\E'$. We define the foliated bundle isomorphism $\psi:f_0^*(E',\E') \times (\RR,\{*\})\to f^*(E',\E')$ over the identity of $M\times \RR$ by setting $\overline \psi(x,v,t):=(x,P_t(v),t)$ where $x\in M, v\in E', t\in \RR, \pi'(v)=f_0(x)$ and $P_t(v)$ is the parallel translation in $E$ of $v$ along $f_t(x)$. The foliated bundles $f_0^*E$ and $f_1^*E$ are isomorphic to the restrictions of the foliated product bundle $f^*E'$ over $M\times [0,1]$ to $M\times \{0\}$ respectively $M\times \{1\}$, and are therefore isomorphic via a foliated bundle map $\varphi$. Then for $(E,\E):=f_0^*(E',\E')$ and $\overline f$ the concatenation of the natural bundle map $f^*E'\to E'$ (having the base map $f$) with $\varphi$ the diagram in Corollary \ref{cor:homotopicbundles}  commutes and by the corollary, up to the isomorphism $\varphi^*$, we have $\overline f_0^*=\overline f_1^*$ on $H_{cv}^*(E,\E)$.

\section{Riemannian foliations and the Molino bundle} \label{sec:Molino}
Let $\F$ be a foliation on a manifold $M$. Recall that by $\Xi(\F)$ we denote the space of differentiable vector fields on $M$ which are tangent to the leaves. A vector field $X$ on $M$ is said to be {\em foliate} if for every $Y\in \Xi(\F)$ the Lie bracket $[X,Y]$ also belongs to $\Xi(\F)$. The flow generated by a foliate field respects the foliation. The set $L(M,\F)$ of foliate fields is the normalizer of $\Xi(\F)$ in the Lie algebra $\Xi(M)$ of vector fields on $M$ and therefore a Lie sub-algebra of $\Xi(M)$. We call the projection of a foliate field $X$ to $TM/T\F$ a {\em transverse} field. The set $l(M,\F)=L(M,\F)/\Xi(\F)$ of transverse fields is also a Lie algebra inheriting the Lie bracket from $L(M,\F)$.

\subsection{Riemannian foliations and the canonical Riemannian basic connection}\label{sec:Riemfol}
Let $\F$ be a foliation and $g$ a metric on $TM/T\F$, called {\em transverse metric}. 
Let $U\subset M$ be a foliation chart domain with projection $p:U\to V\subset \RR^q$ defining $\F$ on $U$ as the fibers. If for every such projection $p$ any two vectors $v,w\in TU/T\F$ with $dp(v)=dp(w)$ have the same length with respect to $g$ we say that $\F$ is a {\em Riemannian foliation}. In this case $V\cong U/\F$ can be endowed with a metric $p_*g$ such that $p:(U,g)\to (V,p_*g)$ becomes a Riemannian submersion. In this sense a Riemannian foliation is locally given by Riemannian submersions. We pull back the Levi-Civita connection from $(V,p_*g)$ to $U$. These connections coincide on overlaps because of the uniqueness of the Levi-Civita connection. This is the canonical Riemannian basic connection. Note that vector fields, that are locally, on a neighborhood $U$ as above, projectable to vector fields on $V$, are parallel along the leaves with respect to this connection.

\subsection{Molino bundle}\label{sec:subMolino}
Let $M$ be a connected $n$-dimensional Riemannian manifold with a Riemannian foliation $\F$ and transverse metric $g$; furthermore, we assume that $M$ is complete with respect to a bundle-like metric that induces $g$ (see Section 3.2 in \cite{Molino}, in particular the last paragraph). We denote the codimension of $\F$ by $q$. Let $\whM$ be the principal $K=\Oq$-bundle over $M$ of orthonormal frames of $\nu\F=TM/T\F$, the {\em Molino bundle} of $(M,\F)$. We denote the natural projection $\whM\to M$ by $\pi$. If $(M,\F)$ is transversely orientable then $\whM$ has two connected $\SOq$-invariant components. A choice of transverse orientation corresponds to a choice of a component. In this case we will by abuse of notation denote this component also by $\whM$ and let $K=\SOq$.  The normal bundle $\nu\F$ is associated to $\whM$, and we denote the connection form on $\whM$ corresponding to the canonical Riemannian basic connection $\nabla$ on $\nu\F$ by $\omega_\F$. We write $H_\pi:=\ker \omega_\F$ for the $K$-invariant horizontal distribution. The manifold $\whM$ carries a natural foliation $\whF$ obtained by horizontally lifting the leaves of $M$. Then $\pi:(\whM,\whF)\to (M,\F)$ is a foliated bundle. Moreover, $\whF$ is respected by the $K$-action, i.e., $K$ maps leaves to leaves.  By construction $\omega$ is a basic form with respect to $\whF$, i.e., $\iota_X\omega_\F=0$ and $\cL_X\omega_\F=0$ for all vector fields $X$ tangential to $\whF$. In particular, we may regard $\omega_\F$ as a map $\omega_\F:\nu \whF\to \soq$. We write $\cH_\pi:=H_\pi/T\whF$ for the {\it transverse horizontal distribution} of $\pi$. Now we lift the transverse metric $g$ on $\nu\F$ to a $K$-invariant metric on the $K$-invariant distribution $\cH_\pi$ of $\nu\whF$. Recall \cite[p.~70, p.~148]{Molino} that the fundamental $1$-form $\theta_\F:\nu\whF\to \RR^q$ is defined by
$$
\theta_\F(X_{\hat x})=\hat x^{-1}(\pi_*(X_{\hat x}))
$$
where $\hat x$ is an orthonormal frame of $\nu_x\F$, understood as the isomorphism $\hat x:\RR^q\to\nu_x\F$ sending the canonical basis to the frame, and $X_{\hat x}\in\nu_{\hat x}\whF= T_{\hat x}\whM/T_{\hat x} \whF$. The fundamental $1$-form $\theta_\F$ is $\whF$-basic by \cite[Lemma 2.1 (i)]{Molino}. By definition, the standard scalar product on $\RR^q$, pulled back with $\theta_\F$ to a degenerate metric on $\nu \whF$, coincides with $\pi^*g$.  We consider the $\whF$-basic, $K$-equivariant map $\omega_\F\oplus\theta_\F:\nu\whF\to\soq\oplus\mR^q$. Pulling back the sum of the standard scalar product on $\RR^q$ and an arbitrary (unique up to a scalar $c$)  biinvariant metric on $\soq$ with $\omega_\F\oplus\theta_\F$ yields an  $K$-invariant $\whF$-transverse metric $\hat g$ on $(\whM,\whF)$ with respect to which $\whF$ is a Riemannian foliation. The projection $\pi$ becomes a Riemannian submersion with respect to the transverse metrics, i.e., $d\pi$ is surjective and the restriction $d\pi:\cH_\pi\to\nu\F$ preserves the metric. We fix the scalar $c$ by requiring that the fibers of $\pi:\whM\to M$ have volume one.

The foliation $\whF$ has a global transverse parallelism, i.e., $\nu\whF$ is parallelizable by transverse fields (we say $\whF$ is TP), see \cite[p.~82, p.~148]{Molino}, with complete representatives in $L(\whM,\whF)$, see \cite[Section 4.1]{GT2010}.

Since $\whF$ is TP the foliation $\overline \whF$ by leaf closures is simple, i.e., $W:=\whM/\overline \whF$ is a manifold and $\overline \whF$ is given as the set of fibers of the locally trivial fibration $\rho:\whM\to W$ (\cite[Proposition 4.1']{Molino}), called the {\em basic fibration}. As the right action of $K$ respects $\whF$ it also respects $\overline \whF$ by continuity and therefore descends to an action on $W$, so that $\rho$ is $K$-equivariant.
The base manifold $W$ can be equipped with a $K$-invariant metric $g_W$ such that $\rho$ becomes a Riemannian submersion with respect to the transverse metric $g_p$ and the metric $g_W$.
Let $\cH_\rho$ the orthogonal complement of $T\overline\whF/T\whF$ in $\nu\whF$ with respect to $\hat g$. We call it the transverse horizontal distribution of $\rho$. Since the transverse metric $\hat g$ is $K$-invariant, so is $\cH_\rho$.

Let us look at the following diagram:
$$
\xymatrix{
  (\whM,\whF,K)\ar[d]^\pi  \ar[r]^\varrho& (W,K)\ar[d] \\
  (M,\F) \ar[r]& M/\oF=W/K.}
$$
Here vertical arrows mean "mod $K$", the horizontal "mod leaf closures". We note that this diagram commutes. In fact, for a leaf closure $N$ of $\F$ the submanifold $N_\whM=\pi^{-1}(N)$ is a $K$-orbit of a leaf closure of $\whF$ therefore descends via $\rho$ to a $K$-orbit $N_W=\rho(N_\whM)$ of $W$. This process can be reversed. Thus, leaf closures of $(M,\F)$ correspond to $K$-orbits of $W$. This is part of a correspondence principle saying that basic data of $\F$ correspond to equivariant data of $(W,K)$. This was exploited in \cite{HT} and in \cite{GT2010}.

\subsection{Molino sheaf and transverse actions}\label{sec:transverseaction}
Let $C(\whM,\whF)$ be the sheaf of local transverse fields that commute with all global transverse fields, called {\em commuting sheaf} of the TP-foliation $\whF$, see \cite[Section 4.4]{Molino}. Its stalk is a Lie algebra $\mfg$ and we write $\umfg:=C(\whM,\whF)$. Any local section of $\umfg$ is the natural lift of a local transverse Killing field of $M$, see  \cite[Proposition 3.4]{Molino}. The {\em Molino sheaf} of $(M,\F)$ is defined as the sheaf on $M$ whose sections are the local transverse Killing fields that naturally lift to local sections of $\umfg$. Thus the Molino sheaf of $(M,\F)$ can be identified with the push-forward $\pi_*C(\whM,\whF)$ of the commuting sheaf. By abuse of notation we denote the Molino sheaf also by $\umfg$. Its stalk is also $\mfg$ which we call the {\em structural Killing algebra} of $(M,\F)$. In Molino's terminology the Molino sheaf is called commuting sheaf of $(M,\F)$, denoted by $C(M,\F)$. A slight difference is that for him, the structural algebra of $(M,\F)$ is the inverse Lie algebra of the stalk of $C(M,\F)$.
For the next definition note that the Molino sheaf is locally constant by \cite[Theorem 5.2]{Molino}.
\begin{dfn}[\cite{Mozgawa}]\label{defn:Killing}
A {\em Killing foliation} is a Riemannian foliation whose Molino sheaf is globally constant.
\end{dfn}
A Riemannian foliation on a simply-connected manifold is therefore automatically a Killing foliation. $(M,\F)$ is a Killing foliation if and only if $C(\whM,\whF)$ is globally constant. In this case the structural Killing algebra $\mfg$ is the center of $l(\whM,\whF)$. Identifying $\mfg$ with the Lie algebra of global sections of the Molino sheaf, $\mfg$ is contained and central in $l(M,\F)$. See also \cite[Theorem 5.2]{Molino}. Therefore $\mfg$ is abelian (but not necessarily the full center of $l(M,\F)$). In order to indicate this we will denote the structural Killing algebra of a Killing foliation by $\mfa$.

\begin{dfn}[{\cite[Sec.~2]{GT2010}}]\label{defn:transverseaction}
A  {\em transverse action} of a finite-dimensional Lie algebra $\mfg$ on a foliated manifold $(M,\F)$ is a Lie algebra homomorphism $\mfg\to l(M,\F),X\mapsto X^*$.
\end{dfn}
The structural Killing algebra $\mfa$ lies in $l(\whM,\whF)$ respectively in $l(M,\F)$ and therefore defines transverse actions on $(\whM,\whF)$ and $(M,\F)$. The projection $\pi$ is $\mfa$-equivariant, i.e.~$d\pi(X^*)=X^*\circ\pi$. The orbits of leaves under the two transverse actions are the respective leaf closures of $\whF$ and $\F$; for more details see \cite[Section 2,Theorem 4.2]{GT2010}.

\section{Basic Thom isomorphism}\label{sec:BasicThom}
Let $\pi:E\to M$ be an oriented fiber bundle with fiber dimension $r$. Let $\pi_*:\Omega^{r+*}_{cv}(E)\to \Omega^*(M)$ be the integration along the fibers, where $\Omega_{cv}$ denotes the de Rham complex of forms with compact vertical support. It is defined as follows. Let $\omega\in\Omega^{r+k}_{cv}(E)$, $p\in M$ and $\xi\in\bigwedge^k(T_pM)$. Then
$$
(\pi_*\omega)_p(\xi):=\int_{E_p}\omega(\,\cdot\, ,\overline\xi),
$$
where $\overline\xi$ is a field of alternating vectors along $E_p$ that map to $\xi$ via $\pi_*$. See for instance \cite[I.7.12]{GHV}; in \cite[p.~61]{BottTu} an equivalent definition is given.

\begin{prop}\label{prop:fiberintegrationprops}
We have
\begin{align*}
&d\pi_*=\pi_*d\\
&\iota_X\pi_*=\pi_*\iota_{\widetilde X}\\
&\cL_X\pi_*=\pi_*\cL_{\widetilde X},
\end{align*}
where $X$ and $\widetilde X$ are $\pi$-related. Moreover, the projection formula
$$
\pi_*((\pi^*\tau)\wedge \omega)=\tau\wedge \pi_*\omega
$$
holds for any form $\tau$ on $M$ and any form $\omega$ on $E$ with compact vertical support.
\end{prop}
\begin{proof}
See \cite[Propositions IX and X, I.7.13]{GHV} and \cite[Proposition 6.15(a)]{BottTu}.
\end{proof}

A foliated bundle $\pi:(E,\E)\to (M,\F)$ is called {\em orientable} if there is a nowhere vanishing $\E$-basic top form. Let our foliated bundle $\pi$ be oriented. By Proposition \ref{prop:fiberintegrationprops} $\pi_*$ respects the basic subcomplexes, and we obtain a chain map $\pi_*:\Omega^{r+*}_{cv}(E,\E)\to \Omega^*(M,\F)$ giving rise to
$$
\pi_*:H^{r+*}_{cv}(E,\E)\to H^*(M,\F).
$$
We note that $\pi^*:\Omega^*(M,\F)\to \Omega^*(E,\E)$ also respects the basic subcomplexes.

In analogy to the general case we have:
\begin{prop}\label{prop:fiberintegrationnaturality}
Integration along the fibers as a homomorphism of basic complexes is natural with respect to foliate bundle maps that are fiberwise diffeomorphisms.
\end{prop}
\begin{proof}
Since integration along the fibers is natural for the ordinary differential complexes by \cite[Corollary I, I.7.12]{GHV}, so is its restriction to the basic subcomplexes.
\end{proof}
Let $X$ be a Riemannian manifold. We call a foliated bundle $\pi:(E,\E)\to (M,\F)$ with typical fiber $X$ and structure group $\Iso(X)$ a {\em Riemannian foliated bundle}. If $X$ is orientable, and if the structure group can be reduced to the group of orientation preserving isometries of $X$ (as structure groups of foliated bundles), the foliated bundle is orientable. If $X$ is a Euclidean space and the transition maps are in the orthogonal group $\mathbf O(X)$ we call $\pi:(E,\E)\to (M,\F)$ a {\em Riemannian foliated vector bundle}. A Riemannian foliated bundle has a natural fiberwise metric. If the bundle in addition has a Riemannian basic connection then the fiberwise metric can be extended to a transverse metric of $(E,\E)$ (similar to the Sasaki metric of the double tangent bundle of a Riemannian manifold) such that $\E$ becomes a Riemannian foliation and $\pi$ a Riemannian submersion for the transverse metrics.

\begin{lemma}\label{lem:Gysinhomog}
Let $(M,\F)$ be a Riemannian foliation with a dense leaf $L$ and let $\pi:(E,\E)\to (M,\F)$ be an oriented Riemannian foliated vector bundle of rank $r$ with a Riemannian basic connection. Then the integration along the fibers $\pi_*:H^{r+*}_{cv}(E,\E)\to H^*(M,\F)$ is an isomorphism.
\end{lemma}
\begin{proof}
First let $L$ be closed and $p\in L$. Then $\Omega(L,\F)=\Omega(\{p\})$ and $\Omega_{cv}(E,\E)=\Omega_c(E_p)^K$, where $K\subset \SO(E_p)$ is the closure of the holonomy group of the foliated bundle at $p$. We have to show that $\int:\Omega_c^{r+*}(E_p)^K\to \Omega^*(\{p\})$ induces an isomorphism in cohomology.  By Poincar\'e duality $\int:H_c^{r+*}(E_p)\to H^*(\{p\})$ is an isomorphism. By averaging over the compact Lie group $K$, we have $H^*(\Omega_c(E_p)^K)=H^*_c(E_p)^K$, which is $\RR$ in dimension $r$, since $K$ preserves the volume form, and otherwise zero. So $\int:H^{r+*}(\Omega_c(E_p)^K)\to H^*(\{p\})$ is an isomorphism.

Now we assume that $L$ is not closed. The next aim is to find a description for $H^*(\oL,\F)$. Let $(\whM,\whF)$ be the Molino bundle and define the fibered product $P=\whM\times_M \SO(E)$ which is a principal bundle over $M$ with structure group $S:={\bf O}(q)\times \SO(r)$. It carries a natural foliation $\F_P:=\whF\times_M\E=\{L_1\times_M L_2\mid L_1\in\whF, L_2\in\E\}$ that is respected by the $S$-action. The projection $\rho:P\to M$ is foliate. The map $\omega_\F\times\omega_E\times \theta_\F\times \theta_E:\nu\E|_M=\nu\F\oplus E\to \soq\times \mathfrak{so}(E) \times\RR^q\times \RR^r$ is foliate and defines a transverse parallelism of $\F_P$; here the $\omega$'s are Riemannian basic connection forms and the $\theta$'s the fundamental $1$-forms. Pulling back an invariant metric on the right side yields a transverse metric for which the foliation is Riemannian (compare with Section \ref{sec:subMolino}). Let $N$ be the leaf closure of a leaf of $\F_P$ over $L$. The restriction $\rho|N:N\to \oL$ is a foliate principal bundle with structure group $K:=S_N=\overline{S_\whL}$. The restricted foliation on $N$ inherits a transverse parallelism from $\F_P$ and is therefore a Lie foliation whose Lie algebra we denote by $\mfg$. According to Fedida there is a covering $p_1:\widetilde N\to N$ whose lifted foliation $\F'$ is given as the fibers of a submersion $p_2:\widetilde N\to G$, where $G$ is the simply-connected Lie group with Lie algebra $\mfg$ (\cite[Theorem 4.2', p.~134]{Molino}). The action by deck transformations $\Gamma$ of $p_1$ respects the foliation and therefore descends along the map $p_2$ to a dense subgroup of $G$ acting from the left. Let us consider the lower line in the following diagram.
$$
\xymatrix{
  G\times \RR^r \ar[d]& p_1^*\rho^*E \ar[r]\ar[l] \ar[d]&\rho^*E  \ar[r] \ar[d]&E \ar[d]^{\pi}\\
  G & \widetilde N \ar[r]^{p_1}\ar[l]_{p_2}& N \ar[r]^\rho& \oL.}
$$
Via $p_1^*$ and $p_2^*$ we have isomorphisms $\Omega(N,\F_P) \cong \Omega(\widetilde N,p_1^*\F_P)^{\Gamma} \cong \Omega(G)^G = \Lambda\mfg^*$ of differential complexes, where $\Lambda\mfg^*$ has the Chevalley-Eilenberg differential $d_\mfg$.
Since $(\Omega(\oL,\F),d)$ is the $K$-basic subcomplex of $\Omega(N,\F_P)$ we obtain
$$
\Omega(\oL,\F)\cong (\Lambda\mfg^*)_{\bas K}.
$$
For the upper line of the diagram note that $\rho^*E$ is a foliated bundle and therefore also $p_1^*\rho^*E$. Note that the action of $\Gamma$ naturally extends to $p_1^*\rho^*E$. In other words, this bundle is $\Gamma$-equivariant. Then it descends along $p_2$ to a $G$-invariant vector bundle over $G$, i.e., to $G\times\RR^r$. This means $\Omega_{cv}(\rho^*E)\cong \Omega_{cv}(G\times\RR^r)^G\cong\Lambda\mfg^*\otimes \Omega_c(\RR^r)$. Passing to the $K$-basic subcomplexes yields
$$
\Omega_{cv}(E,\E)\cong (\Lambda\mfg^* \otimes \Omega_c(\RR^r))_{\bas K}.
$$
Here $\bas K$ means the $\mfk$-horizontal, $K$-invariant part. We want to see that the integration along the fibers $\pi_*:\Omega_{cv}(E,\E)\cong (\Lambda\mfg^* \otimes \Omega_c(\RR^r))_{\bas K}\to (\Lambda\mfg^*)_{\bas K}\cong \Omega(M,\F)$ induces an isomorphism in cohomology. Define the two $K$-differential graded algebras $A=\Lambda\mfg^*$ and $B=\Omega_c(\RR^r)$ (for the definition see \cite[Definition 3.12]{kambertondeur} or \cite[Definition 2.3.1]{GS}; compare with Definition \ref{defn:g*}). Let $\xi_a, a=1,\ldots,\dim\mfk$ be a basis of $\mfk$, the Lie algebra of $K$. Since $K$ is compact, $\mfg$ has an ${\rm Ad}(K)$-invariant subpace $\mathfrak p$ complementary to $\mfk$. We extend the basis of $\mfk$ by $\xi_b\in \mathfrak p, b=\dim\mfk+1,\ldots,\dim\mfg$ to a basis of $\mfg$. Let $\theta^1,\ldots,\theta^{\dim\mfg}$ be the dual basis of $\xi_1,\ldots,\xi_{\dim\mfg}$. From now on the variable $a$ runs in the range $1,\ldots,\mfk$, $b$ in the range $\dim\mfk+1,\ldots,\dim\mfg$ and $c$ in the range $1,\ldots, \dim\mfg$. The set of $\theta^a$ is dual to the set of $\xi_a$ and their span in $A$ is $K$-invariant; in other words $A$ satisfies condition $(C)$, see \cite[Definition 2.3.4]{GS} (this is equivalent to the existence of a connection in the sense of \cite[Definition 3.4]{Meinrenken}). We will write $\iota_a$ and $\cL_a$ for the derivations $\iota_{\xi_a}$ and $\cL_{\xi_a}$. For two $K$-differential graded algebras $A$ and $B$ with $A$ satifisfying condition $(C)$ there is a $K$-equivariant algebra automorphism $\phi$ of $A\otimes B$, the {\em Mathai-Quillen isomorphism}, mapping $(A\otimes B)_{\rm hor}$ to $A_{\rm hor}\otimes B$ (see \cite[(4.12)]{GS}), where $\rm hor$ denotes the $\mfk$-horizontal part, and therefore $(A\otimes B)_{\bas K}$ to $(A_{\rm hor}\otimes B)^K$. The Mathai-Quillen isomorphism $\phi:A\otimes B\to A\otimes B$ is defined by $\phi:=\exp\gamma=1+\gamma+\frac12\gamma^2+ \ldots+ \frac{1}{(\dim\mfk)!}\gamma^{\dim\mfk}$ where $\gamma:=\sum_a\theta^a\otimes \iota_a\in {\rm End}(A\otimes B)$. In our case we obtain the isomorphism
$$
\phi: {(\Lambda \mfg^*\otimes \Omega_c(\RR^r))}_{\bas K}\to ({(\Lambda \mfg^*)}_{\rm hor}\otimes \Omega_c(\RR^r))^K.
$$
The differential $d=d_\mfg\otimes 1 + 1 \otimes d_{\RR^r}$  translates via $\phi$ into
\begin{align*}
D:=&\phi d\phi^{-1}=d_\mfg\otimes 1 + 1 \otimes d_{\RR^r} - \sum_a\mu^a\otimes\iota_a+\sum_a\theta^a\otimes \cL_a\\
=&\sum_a\theta^a(\cL_a\otimes 1 + 1\otimes \cL_a) + 1 \otimes d_{\RR^r} - \sum_a\mu^a\otimes \iota_a +\sum_b\theta^b\otimes \cL_b\\
=&1 \otimes d_{\RR^r} - \sum_a\mu^a\otimes \iota_a + \sum_b\theta^b\otimes \cL_b
\end{align*}
The second equality is \cite[Equation (4.11)]{GS}, the third follows from $d_\mfg=\sum_c \theta^c \cL_c$, the fourth since $\cL_a\otimes 1 + 1\otimes \cL_a$ is zero on the invariant forms; here $\mu^a$ are the curvature forms (\cite[Equation (3.10)]{GS}).
By passing to $K$-invariant forms we obtain the top horizontal line in the following diagram of differential complexes.
$$
\xymatrix{
  \Omega_{cv}(E,\E)\ar[d]^{\pi_*} \ar[r]^\cong & {(\Lambda \mfg^*\otimes \Omega_c(\RR^r))}_{\bas K} \ar[d]^{s}\ar[r]_{\phi}^\cong & {((\Lambda \mfg^*)}_{\rm hor}\otimes \Omega_c(\RR^r))^K \ar[d]^{s}\\
  \Omega_{cv}(M,\F)\ar[r]^\cong & {(\Lambda \mfg^*)}_{\bas K}\ar[r]^= &{(\Lambda \mfg^*)}_{\bas K}.}
$$
We have to identify the middle and right vertical maps and justify the destination space on the right. The map $s:A\otimes B\to A\otimes \RR=A$ that integrates the second factor $B=\Omega_c(\RR^r)$ is a morphism of $K$-differential graded algebras (note that all operators $\iota_a,\cL_a, d$ are zero on $\RR$) and descends to a map of $K$-basic subcomplexes that makes the first square in the diagram commutative. For the right vertical map we observe $s\gamma=0=\gamma s$ thus $s$ commutes with $\gamma$ and consequently with $\phi$. So the right vertical map $((\Lambda \mfg^*)_{\rm hor} \otimes \Omega_c(\RR^q))^K\to (\Lambda \mfg^*)_{{\rm bas}\, K}$ is the restriction of $s$.

We will now show that the map $s$ on the right induces an isomorphism in cohomology. Let us consider the bigraded complex $C^{s,t}:=((\Lambda^s\mfg^*)_{\rm hor} \otimes \Omega^t_c(\RR^r))^K$. The differential $D$ is then the sum of maps of degree $(0,1)$, $(2,-1)$ and $(1,0)$. With $F_p:=\bigoplus_{s\geq p,t\in\NN}C^{s,t}$ we have a filtration $F_0\supset F_1 \supset \ldots \supset F^{\dim\mfg+1}=0$ that is respected by $D$. The $E_1$-term of the associated spectral sequence is $E_1=H_{1\otimes d_{\RR^r}}(C^{\bullet,\bullet})=((\Lambda^\bullet\mfg^*)_{\rm hor}\otimes H_c^\bullet(\RR^r))^K$. Similarly we filter the complex $(\Lambda\mfg^*)_{\bas K}$ by $\widetilde F_p:=\bigoplus_{s\geq p}(\Lambda^s\mfg)_{\bas K}$. Note that both associated spectral sequences converge because they are of finite length. The integration
$s:((\Lambda \mfg^*)_{\rm hor} \otimes \Omega_c(\RR^r))^K\to (\Lambda \mfg^*)_{\bas K}$ of the second factor is a cochain map and preserves the filtrations. Since it induces an isomorphism at the $E_1$-stage, it also induces isomorphisms at the $E_l$-stage for all $l\geq 1$ (see \cite[Theorem 3.4]{McC}) and therefore yields an isomorphism
$$
s:H(((\Lambda \mfg^*)_{\rm hor}\otimes \Omega_c(\RR^r))^K,D)\to H((\Lambda \mfg^*)_{\bas K},d_\mfg)
$$
So $\pi_*:H^{r+*}_{cv}(E,\E)\to H^*(M,\F)$ is an isomorphism.
\end{proof}
\begin{dfn}\label{def:finitetype}
Let $(M,\F)$ be a foliated manifold. A cover $\cU=\{U_\alpha\}$ of $M$ by open saturated sets $U_\alpha$ is called a {\em good saturated cover} if all nonempty finite intersections $V=U_{\alpha_0}\cap \cdots \cap U_{\alpha_p}$ are $\F$-contractible, meaning that there is a leaf closure $\oL\subset V$ and a basic homotopy $H:V\times [0,1]\to V$ with $H_0={\rm id}_V$ and $H_1(V)\subset \oL$.
We call $(M,\F)$ of {\em finite type} if it has a finite good saturated cover.
\end{dfn}
\begin{prop}\label{prop:finitetype}
 A Riemannian foliation $\F$ on a manifold $M$ that is complete with respect to an adapted bundle-like metric has a good saturated cover. If $M/\oF$ is compact, $(M,\F)$ is of finite type. Moreover, for any good saturated cover there is a basic partition of unity subordinate to it.
\end{prop}
\begin{proof}
See appendix.
\end{proof}
\begin{thm}\label{thm:thomiso}
Let $(M,\F)$ be a Riemannian foliation of finite type and let $\pi:(E,\E)\to (M,\F)$ be an oriented Riemannian foliated vector bundle of rank $r$ with a Riemannian basic connection. Then the integration along the fiber
$$
\pi_*:H^{r+*}_{cv}(E,\E)\to H^*(M,\F).
$$
is an isomorphism.
\end{thm}
\begin{proof}
The general structure of the proof is as that of \cite[Prop 6.17]{BottTu}. We will start with a number of observations. Let $U$ be an open saturated subset of $M$ that is $\F$-contractible to a leaf closure $\oL$. We claim that $\pi_*:H^{r+*}_{cv}(E|_U,\E)\to H^*(U,\F)$ is an isomorphism. In fact, let $H:U\times [0,1]\to U$ be a basic homotopy with $H_0={\rm id}_U$ and $H_1(U)\subset \oL$. We define $\overline H:E|_U\times [0,1]\to E|_U$ such that $\overline H_t(v)$ for $v\in E_x,x\in U$ is the parallel translation of $v$ along the curve $H_t(x)$. Then the following diagram commutes by the naturality of $\pi_*$, Proposition \ref{prop:fiberintegrationnaturality}.
$$
\xymatrix{
  H^{r+*}_{cv}(E|U,\E)\ar[d]^{\pi_*} & H^{r+*}_{cv}(E|{\oL},\E)  \ar[d]^{\pi_*}\ar[l]_{\overline H^*_1} \\
  H^*(U,\F) & H^*(\oL,\F)\ar[l]^{H^*_1} .}
$$
The lower map is an isomorphism because of Corollary \ref{cor:homotopicmaps}, the upper one because of the discussion succeeding Corollary \ref{cor:homotopicbundles}, and the right map because of Lemma \ref{lem:Gysinhomog}. This proves the claim.

Let $U$ and $V$ be saturated open subsets of $(M,\F)$. We want to see that if $\pi_*$ is an isomorphism over $U, V$ and $U\cap V$ then also over $U\cup V$. Using a basic partition of unity subordinate to $\{U,V\}$, which exists by Proposition \ref{prop:finitetype}, we see that the sequence
$$
0\to \Omega^*_{cv}(E|({U\cup V}),\E)\longrightarrow \Omega^*_{cv}(E|{U},\E) \oplus  \Omega^*_{cv}(E|{V},\E) \longrightarrow \Omega^*_{cv}(E|({U\cap V}),\E)\to 0
$$
is exact. So we have the Mayer-Vietoris sequences
\[\tiny{
\xymatrix@R=20pt@C=10pt{
H^*_{cv}(E|({U\cup V}),\E) \ar[d]^{\pi_*} \ar[r] & H^*_{cv}(E|{U},\E) \oplus  H^*_{cv}(E|{V},\E) \ar[d]^{\pi_*} \ar[r] & H^*_{cv}(E|({U\cap V}),\E) \ar[d]^{\pi_*} \ar[r] & H^{*+1}(E|({U\cup V}),\E) \ar[d]^{\pi_*} \\
 H^*(U\cup V,F)\ar[r]  &  H^*(U,\F) \oplus  H^*(V,\F) \ar[r] &  H^*(U\cap V,\F) \ar[r] &  H^{*+1}(U\cup V,\F)   } }
\]
The commutativity is clear for the first two squares and is proven as in the classical case for the third. So if $\pi_*$ is an isomorphism over $U,V$ and $U\cap V$ then by the Five-Lemma also over $U\cup V$. Now we can prove the theorem by induction on the cardinality $p$ over a good saturated cover. For $p=1$ this is true by the discussion at the beginning of the proof. Now suppose $\pi_*$ is an isomorphism for any Riemannian foliation having a good saturated cover with at most $p$ open sets. Consider a Riemannian foliation $(M,\F)$ having a good saturated cover $\{U_0,\ldots U_{p+1}\}$ with $p+1$ open sets. Now $(U_0\cup \ldots \cup U_{p-1})\cap U_p$ has a good saturated cover with cardinality $p$, namely $\{U_0\cap U_p, \ldots, U_{p-1}\cap U_p\}$. By hypothesis $\pi_*$ is an isomorphism over $(U_0\cup\ldots \cup U_{p-1})\cap U_p$ and over $U_p$, and therefore, by the previous discussion, over the union $U_0\cup \ldots \cup U_p$.
\end{proof}
\begin{remark}
From a more formal point of view in the proof we have actually considered $\pi_*$ as a homomorphism
$$
\pi_*:C^{p,q+r}(\pi^{-1}(\cU),\Omega(\,\cdot\,,\E))\to C^{p,q+r}(\cU,\Omega(\,\cdot\,,\F))
$$
between Cech-basic de Rham complexes.
\end{remark}

In the situation of the theorem the inverse map $\tau:H^*(M,\F)\to H^{r+*}_{cv}(E,\E)$ of the integration along the fiber $\pi_*:H^{r+*}_{cv}(E,\E)\to H^*(M,\F)$ is called {\em basic Thom isomorphism}. The image $\Phi$ in $H^{r}_{cv}(E,\E)$ of $1\in H^0(M,\F)$ is called {\em basic Thom class} of $(E,\E)$. Because $\pi_*\Phi=1$, by the projection formula in Proposition \ref{prop:fiberintegrationprops}
$$
\pi_*(\pi^*\omega\wedge \Phi)=\omega\wedge\pi_*\Phi=\omega
$$
for all $\omega\in H^*(M,\F)$.
So the Thom isomorphism is given just as in the classical case by
$$
\tau(\ )=\pi^*(\ )\wedge \Phi.
$$
Moreover, as in \cite[Proposition 6.18]{BottTu} one can show that the basic Thom class is uniquely characterized as the basic cohomology class in $H^{r}_{cv}(E,\E)$ that maps to the generator of $H_c(E_p)$ for each fiber $E_p$.

Let $(M,\F)$ be a Riemannian foliation and $\nu\F$ be its normal bundle. It carries a natural foliation $\whF$. Indeed, the leaf $\whL_v$ through a vector $v\in\nu_x\F$ is the submanifold consisting of all normal vectors that are obtained by holonomy translation of $v$ along the leaf $L_x$. By construction the projection $\pi:(\nu\F,\whF)\to (M,\F)$ is foliate and the normal bundle is a foliated bundle.
The transverse metric of $\F$ makes $(\nu\F,\whF)$ into a Riemannian foliated bundle over $M$. Let $N$ be a closed stratum of the $\mfa$-stratification introduced in \cite[Section 4.2]{GT2010}, e.g.~a connected component of the union of all closed leaves, and let $r$ be its codimension in $M$. The Riemannian foliation $\whF$ restricts a Riemannian foliation on $\nu N$ which we also denote by $\whF$. We observe that the canonical Riemannian basic connection of $\F$ preserves the normal bundle of $N$.
\begin{cor}\label{cor:basicThomstratum}
Let $\F$ be a Riemannian foliation on a complete manifold $M$. If the foliated bundle $(\nu N,\whF|\nu N)$ is oriented (e.g.~if $\F$ is Killing) then the integration along the fibers $\pi_*:H_{cv}^{r+*}(\nu N,\whF|\nu N)\to H^{*}(N,\F)$ is an isomorphism.
\end{cor}
\begin{proof}
It remains to show that $\nu N$ is orientable if $\F$ is a Killing foliation. In this case the transverse action of the structural Killing algebra $\mfa$ is isometric, see \cite[Section 2.1]{GT2010} and therfore restricts on each normal space $\nu_pN$ to an action by skew-symmetric automorphisms. Since $\mfa$ is abelian, $\nu_pN$ decomposes into a direct sum of its two-dimensional weight spaces. A choice of a generic element of $\mfa$ orients these weight spaces simultaneously (an appropriate normalization of this skew-symmetric automorphism even defines an $\mfa$-invariant complex structure for these weight spaces). This provides an orientation of $\nu N$ as a foliated bundle.
\end{proof}
We will call the concatenation
$$
i_*^N:H^*(N,\F)\stackrel{\tau}{\to} H_{cv}^{r+*}(\nu N,\whF|\nu N) \to H^{r+*}(M,\F)
$$
of the basic Thom isomorphism of the foliated bundle $(\nu N,\whF|\nu N)$ with the homomorphism induced by the inclusion of the normal bundle, seen as a saturated tubular neighborhood of $N$, into $M$ the {\em basic Thom homomorphism} of $(N,\F)$ in $(M,\F)$.

\section{Transverse integration and basic Poincar\'e duality}\label{sec:PoincareDuality}
We briefly recall Poincar\'e duality for the basic cohomology ring of a Riemannian foliation according to \cite{Sergiescu}. We will first assume that $\F$ is a transversely oriented Killing foliation of codimension $q$ on a complete manifold $M$. Let $\pi:(\whM,\whF)\to (M,\F)$ be the $\SOq$-Molino bundle with structural Killing algebra $\mfa$. For $X\in\mfa$ the corresponding transverse field of the transverse action by the structural Killing algebra is denoted by $X^*\in l(M,\F)$. Moreover, we set $k:=\dim\mfa$.

Sergiescu defines an operator $\int_\F$ on $\Omega_c(M,\F)$, the complex of basic forms of $(M,\F)$ whose support projected to $M/\oF$ is compact. We fix a nontrivial linear $k$-form $\nu_\mfa\in\Lambda^k(\mfa^*)$. Let $X_1,\ldots,X_k$ be basis of $\mfa$ such that $\nu_\mfa(\mathbf X)=1$ for $\mathbf X=X_1\wedge \ldots \wedge X_k\in \Lambda^k(\mfa)$ and define the operator $\iota_{\mathbf X}:=\iota_{X_k^*}\circ \cdots \circ \iota_{X_1^*}:\Omega^{*+k}(\whM,\whF)\to\Omega^*(\whM,\whF)$. Now $\nu_\mfa$ determines a basic $k$-form $\nu_\rho$ of $(\whM,\whF)$ by requiring $\iota_{\mathbf X}\nu_\rho=1$ and $\iota_Y\nu_\rho=0$ for every $Y\in \cH_\rho$ where $\cH_\rho$ is the transverse horizontal distribution of $\rho$ (see Section \ref{sec:Molino}).
\begin{lemma}\label{lem:rho}
For $\omega\in\Omega(\whM,\whF)$ the form $\iota_{\mathbf X}\omega$ is basic with respect to $\rho$. Moreover, $d\iota_{\mathbf X}=(-1)^{\dim\mfa}\iota_{\mathbf X}d$.
\end{lemma}
\begin{proof}
The form $\eta:=\iota_{\mathbf X}\omega$ is $\whF$-basic, see \cite[p.~39]{Molino}.  We show that $\eta$ is also basic with respect to the transverse $\mfa$-action. The form $\eta$ is $\mfa$-horizontal since $\iota_{X^*}\iota_{\mathbf X}=0$ for all $X\in\mfa$, and it is $\mfa$-invariant as an $\whF$-basic form (see \cite[Lemma 3.15]{GT2010}). Now recall that the fibers of $\rho$ are the leaf closures of $\whF$ which are orbits of leaves under $\mfa$. Therefore $\eta$ is basic with respect to all vector fields tangential to $\overline\whF$, hence basic with respect to $\rho$.
The second statement in the lemma follows from the fact that $d\iota_{X^*}\beta+\iota_{X^*}d\beta=\cL_{X^*}\beta=0$ for all $\beta\in\Omega(\whM,\whF)$, again because every $\whF$-basic form is $\mfa$-invariant .
\end{proof}
We can therefore define $\rho_\#:\Omega^{*+k}(\whM,\whF)\to \Omega^*(W)$ by $\rho_\#(\omega)=(\rho^*)^{-1}\iota_{\mathbf X}\omega$ and the transverse integration $\int_\whF:\Omega_c(\whM,\whF)\to \RR$ of $\whF$ by
$$
\int_\whF:=\int_W\circ\ \rho_\#.
$$
Note that this is more or less the map $I$ defined in the proof of \cite[Lemma 2.1]{Sergiescu}. It can only be nonzero on top basic forms. Because of Lemma \ref{lem:rho} and the Theorem of Stokes for $\int_W$ we have $\int_\whF\circ\,d=0$. Since $\F$ is a transversely oriented Killing foliation, $W$ is oriented. So $\int_\whF\nu_\rho\wedge\rho^*\eta=\int_W\eta\neq 0$ for some $\eta\in\Omega^{\dim W}_c(W)$, thus defining a surjection $\int_\whF:H^{\hat q}(\whM,\whF)\to\RR$ where $\hat q:=\codim\whF$. Sergiescu shows that this is in fact an isomorphism (\cite[Proposition 2.4]{Sergiescu}).

Let $\pi_*:H^{*+\dim\soq}(\whM,\whF)\to H^*(M,\F)$ be the integration along the fibers defined in Section \ref{sec:BasicThom}. It is an isomorphism for $*=q$ with inverse $\omega\mapsto \pi^*\omega\wedge\nu_{\pi}$ (which can be verified with the projection formula in Proposition \ref{prop:fiberintegrationprops}), where $\nu_{\pi}$ is the volume element on $\soq$ with respect to the biinvariant inner product chosen in Section \ref{sec:Molino} composed with the connection form $\omega_\F$ (see the proof of \cite[Lemma 2.6]{Sergiescu}). Next we define the transverse integration operator $\int_\F:H^q_c(M,\F) \to\RR$ by
$$
\int_\whF=\int_\F\circ\ \pi_*,
$$
see \cite[Lemma 2.6]{Sergiescu}. Note that $\int_\F\omega=\int_\whF\pi^*\omega\wedge\nu_{\pi}$. It follows
\begin{prop}\label{prop:integration}
Tranverse integration
$$
\int_\F:H^q_c(M,\F)\to\RR
$$
is an isomorphism
\end{prop}
So far we have assumed that $\F$ is a transversely oriented Killing foliation for the definition of the transverse integration operators $\int_\whF$ and $\int_\F$. These can be defined more generally for a Riemannian foliation $\F$ whose sheaf $\Lambda^{\dim\mfg}\umfg^*$ is constant, where $\umfg$ is the commuting sheaf; such an $\F$ is called {\em taut}. Under this condition $\nu_\mfg$ can still be extended to $\Omega(\whM,\whF)$, allowing first to define $\rho_*$ even though $\umfg$ might not be constant, and then $\int_\whF$ and $\int_\F$. For this more general class of Riemannian foliations Proposition \ref{prop:integration} and basic Poincar\'e duality as follows hold:

\begin{thm}[Sergiescu] For a transversely oriented taut Riemannian foliation on a complete manifold the pairings
\begin{align*}
H^*(M,\F)\otimes H^{q-*}_c(M,\F)\to\RR \\
H^*_c(M,\F)\otimes H^{q-*}(M,\F)\to\RR
\end{align*}
defined by $([\omega],[\eta])\mapsto \int_\F \omega\wedge\eta$ are non-degenerate.
\end{thm}

Sergiescu shows basic Poincar\'e duality with twisted basic cohomology for arbitrary Riemannian foliations. Here we will not be concerned with this case.

We remark that in analogy to the standard setting (see e.g.~\cite[p.~50 ff.]{BottTu}) one can now define the basic Poincar\'e dual of a closed, saturated, transversely oriented submanifold of $(M,\F)$.

Let $C$ be the union of closed leaves. We fix a connected component $C_i$ and let $i_*:H^*(C_i,\F)\to H^{r_i+*}(M,\F)$ be the basic Thom homomorphism of $C_i$. Choose $x_i\in C$ with $L_{x_i}$ regular, i.e.~without holonomy, in $(C_i,\F)$, $\hat x_i\in \whM$ with $\pi(\hat x_i)=x_i$ and define $\bar x_i:=\rho(\hat x_i)$. Let
$$
c_i:=s\int_{\SOq/\SOq_{\overline\whL}}\iota_{\mathbf X}\nu_\pi=s\int_{\SOq/\SOq_{\bar x_i}}\rho_\#\nu_\pi
$$
where $\whL$ is the leaf of $(\whM,\whF)$ through $\hat x$ and $s=(-1)^{\dim\soq+\codim C_i+\dim C_i/\F}$. Recall from Section \ref{sec:subMolino} that $\rho:\whW\to W$ is $\SOq$-equivariant.
\begin{lemma}\label{lem:restrictionintegral}
We have $\int_\F\circ\, i_*=c_i\int_{C_i/\F}$, where the latter denotes integration over the orbifold $C_i/\F$.
\end{lemma}
\begin{proof}
Let $\cH_\pi\subset\nu\whF$ be the transverse horizontal distribution of the transversely Riemannian submersion $\pi:\whM\to M$ and $\cH_\rho\subset\nu\whF$ the transverse horizontal distribution for $\rho:\whM\to W$. We define $N_\whM=\pi^{-1}(N)$ and $N_W=\rho(N_\whM)$ for any saturated set $N\subset M$. For $N:=C_i$ we have
$$
(TN/T\F)^{\cH_\pi}_{\hat x}\oplus \hat x\cdot  \soq= T_{\hat x}N_\whM/T_{\hat x} \whF= (TN_W)^{\cH_\rho}_{\hat x}\oplus \mfa\cdot\hat x,
$$
by \cite[Proposition 4.7]{GT2010}, where the direct sums are orthogonal with respect to the transverse metric $\hat g$ and the superscripts denote the respective horizontal lifts, e.g.~$(TN/T\F)^{\cH_1}$ is the geometric realization of $\pi^*(TN/T\F)$ as a subbundle of $\cH_1$. As $\pi$ and $\rho$ are transversely Riemannian submersions, $\pi_*:\cH_\pi\to \nu \F$ and $\rho_*:\cH_\rho\to TW$ are pointwise isometric, and we have
$$
(\pi^*\nu N)_{\hat x}\cong(\nu N)^{\cH_\pi}_{\hat x}= \nu_{\hat x} N_\whM= (\nu N_W)^{\cH_\rho}_{\hat x}\cong (\rho^*\nu N_W)_{\hat x}.
$$
Here $\nu N=(TN/T\F)^{\perp_g}, \nu N_\whM= (TN_\whM/T\whF)^{\perp_{\hat g}}$ and $\nu N_W=(TN_W)^{\perp_{g_W}}$. Let $U$ be a small saturated tubular neighborhood of $N$ and let $p:U\to N$, $\hat p:U_{\whM}\to N_\whM$ and $\bar p:U_W\to N_W$ the geodesic projection maps with $\pi\circ\hat p=p\circ\pi$ and $\rho\circ\hat p=\bar p\circ\rho$.
Let $\omega$ be a basic form of $(U,\F)$. Let $\Phi$ be the basic Thom class of the foliated bundle $(\nu N,\whF)\cong (U,\F)$ and $\Phi_W$ the Thom class of $\nu N_W$. Since a (basic) Thom class is represented by any closed (basic) form that restricts to the generator of the cohomology of a fiber (see Section \ref{sec:BasicThom}), $\pi^*\Phi$ and $\rho^*\Phi_W$ are the basic Thom class of $\nu N_{\whM}$, so
$$
\pi^*\Phi=\rho^*\Phi_W
$$ 
in cohomology. Since $\hat p$ is $\SOq$- and $\mfa$-equivariant one has $\hat p^*\nu_\pi|_{N_\whM}=\nu_\pi$ and $\hat p^*(\nu_\rho|_{N_\whM})=\nu_\rho$. We define  $\eta=(\rho^*)^{-1}\pi^*\omega\wedge\rho_\#\nu_\pi|_{N_\whM}$; note that $\pi^*\omega$ is $\mfa$-horizontal by the first highlighted equation above and $\mfa$-invariant as an $\whF$-basic form, therefore $\rho$-basic. Let $k=\dim\mfa=\deg \nu_\rho, l=\dim\soq=\deg\nu_\pi, m=\codim C_i, w=\dim C_i/\F$. Then $\pi^*\omega\wedge\nu_\pi|_{N_\whM}=(-1)^l\rho^*\eta\wedge\nu_\rho|_{N_\whM}$, because $\nu_\pi=\nu_\rho\wedge\iota_{\mathbf X}\nu_\pi$. Furthermore
\begin{align*}
\pi^*(p^*\omega\wedge\Phi)\wedge\nu_\pi =& (-1)^{l+m}\pi^*p^*\omega\wedge\nu_\pi \wedge\pi^*\Phi\\
=& (-1)^{l+m}\hat p^*(\pi^*\omega\wedge\nu_\pi|_{N_\whM})\wedge \pi^*\Phi\\
=& (-1)^{m}\hat p^*(\rho^*\eta\wedge\nu_\rho|_{N_\whM})\wedge \rho^*\Phi_W\\ 
=& (-1)^{m} \rho^*\bar p^*\eta\wedge \nu_\rho\wedge \rho^*\Phi_W\\
=& (-1)^{m+w+k}\nu_\rho\wedge \rho^*(\bar p^*\eta\wedge \Phi_W). \\
\end{align*}
Now
\begin{align*}
\int_\F i_*\omega =& \int_\F p^*\omega\wedge\Phi=\int_\whF \pi^*(\bar p^*\omega\wedge\Phi)\wedge\nu_\pi\\
=&(-1)^{m+w+k}\int_\whF \nu_\rho\wedge\rho^*(\bar p^*\eta\wedge \Phi_W) \\
=&(-1)^{m+w+k}\int_W \bar p^*\eta\wedge \Phi_W\\
=&(-1)^{m+w+k} \int_{(C_i)_W} \bar p^*\eta\\
=&(-1)^{m+w+k}\int_{(C_i)_W} (\rho^*)^{-1}\pi^*p^*\omega\wedge\rho_\#\nu_\pi\\
=&c_i \int_{(C_i)_W/\SOq} (\rho^*)^{-1}\pi^*p^*\omega\\
=&c_i \int_{C_i/\F}\omega.\\
\end{align*}
\end{proof}

\section{Equivariant basic cohomology}\label{sec:EquivariantBasicCohomology}
In this section we will recall the notion of equivariant basic cohomology from \cite{GT2010}.
\begin{dfn}\label{defn:g*}
Let $\mfg$ be a finite-dimensional Lie algebra and $A=\bigoplus A_k$ a $\ZZ$-graded algebra. We call $A$ a
$\mfg$-differential graded algebra ($\mfg$-dga) if there is a derivation $d:A\to A$ of degree $1$ and  derivations $\iota_X:A\to A$ of degree $-1$ and $\cL_X:A\to A$ of degree $0$ for all $X\in \mfg$ (where $\iota_X, \cL_X$ linearly depend on $X$) such that:
\begin{enumerate}
\item $d^2=0$
\item $[\cL_X,L_Y]=\cL_{[X,Y]}$
\item $[\cL_X,\iota_Y]=\iota_{[X,Y]}$
\item $[d,\cL_X]=0$
\item $\iota_X\iota_Y+\iota_Y\iota_X=0$
\item $d\iota_X+\iota_Xd=\cL_X$.
\end{enumerate}
\end{dfn}

\begin{ex} \label{ex:gstaraction}
An  infinitesimal action of a finite-dimensional  Lie algebra $\mfg$ on a manifold $M$, i.e.~a Lie algebra homomorphism $\mfg\to\Xi(M)$, induces a $\mfg$-dga structure on the de Rham complex $\Omega^*(M)$.
\end{ex}

We want to define the equivariant cohomology of an arbitrary $\mfg$-dga $A$. First define the {\em Cartan complex}
\begin{equation*}
C_\mfg(A):=(S(\mfg^*)\otimes A)^\mfg.
\end{equation*}
Here the superscript denotes the subspace of $\mfg$-invariant elements, i.e., those $\omega\in S(\mfg^*)\otimes A$ for which $L_X\omega=0$ for all $X\in \mfg$. The differential $d_\mfg$ of the {\em Cartan complex} $C_\mfg(A)$ is defined by
$$
(d_\mfg \omega)(X)=d(\omega(X))-\iota_X(\omega(X)),
$$
where we consider an element in $C_\mfg(A)$ as a $\mfg$-equivariant polynomial map $\mfg\to A$. Choosing a basis $\{X_i\}_{i=1,\ldots,r}$ of $\mfg$ with dual basis $\{u_i\}_{i=1,\ldots,r}$ of $\mfg^*$ we have
$$
d_\mfg\omega:=d\omega-\sum_{i=1}^r\iota_{X_i}(\omega)u_i.
$$
Since $d_\mfg\circ d_\mfg=0$ on $C_\mfg(A)$ we can define the equivariant cohomology of the $\mfg$-dga $A$ by
\begin{equation}\label{eq:cartandef}
H_\mfg^*(A):=H^*(C_\mfg(A),d_\mfg).
\end{equation}
Note that there is a natural $S(\mfg^*)^\mfg$-algebra structure on $H^*_\mfg(A)$.

A graded algebra homomorphism $f:A\to B$ between $\mfg$-dgas intertwining  $d, \iota_X$ and $\cL_X$ for all $X\in \mfg$ is called a $\mfg$-dga-homomorphism. Such a homomorphism induces a chain map $f_*:(C_\mfg(A),d_\mfg)\to (C_\mfg(B),d_\mfg)$ between the corresponding Cartan complexes and therefore an algebra homomorphism $f_*:H_\mfg(A)\to H_\mfg(B)$ between the corresponding equivariant cohomology algebras, which is moreover an $S(\mfg^*)^\mfg$-module homomorphism.

We will now apply this concept to transverse actions on foliations. Recall from Definition \ref{defn:transverseaction} that a transverse action of a Lie algebra $\mfg$ on a  foliated manifold $(M,\F)$ is a Lie algebra homomorphism $\mfg\to l(M,\F)$. For $X\in\mfg$ we will denote the corresponding transverse field by $X^*\in l(M,\F)$. Let $\tilde X\in L(M,\F)$ be a foliate field that represents $X^*$. The derivations $\iota_X:=\iota_{\tilde X}$ and $\cL_X:=\cL_{\tilde X}$ as operators on $\Omega^*(M,\F)$ do not depend on the choice of representative $\tilde X$. Together with the restriction of the differential of $\Omega^*(M)$ we have the following.

\begin{prop}[{\cite[Prop.~3.12]{GT2010}}]\label{prop:omegaMFgstar}
A transverse action of a finite-di\-mensio\-nal  Lie algebra $\mfg$ on a foliated manifold $(M,\F)$ induces the structure of a $\mfg$-dga on $\Omega^*(M,\F)$.
\end{prop}

This proposition enables us to apply the general construction of equivariant cohomology of a $\mfg$-dga as defined in \eqref{eq:cartandef} to $\Omega(M,\F)$. We will write $\Omega_\mfg^*(M,\F)$ for $C_\mfg^*(\Omega(M,\F))$ and call its elements {\em equivariant basic forms}.

\begin{dfn}[{\cite[Def.~3.13]{GT2010}}] For a foliated manifold $(M,\F)$ with a transverse action of a finite-dimensional  Lie algebra $\mfg$ we define the {\em equivariant basic cohomology} of the $\mfg$-action on $(M,\F)$ by
 $$
H_\mfg(M,\F):=H(\Omega_\mfg(M,\F)),d_\mfg)=H((S(\mfg^*)\otimes \Omega(M,\F))^\mfg),d_\mfg).
$$
\end{dfn}

Now let $\pi:(E,\E)\to (M,\F)$ be a foliated vector bundle and assume there are transverse actions of a Lie algebra $\mfg$ on $(E,\E)$ and $(M,\F)$ such that $\pi$ is $\mfg$-equivariant, i.e., $d\pi( X^*)=X^*\circ \pi$ for all $X\in\mfg$, see Section \ref{sec:transverseaction}. By Proposition \ref{prop:fiberintegrationprops} we obtain a $\mfg$-dga homomorphism $\pi_*:\Omega^{r+*}_{cv}(E,\E)\to \Omega^*(M,\F)$ of degree $-r$, inducing an $S(\mfg^*)^\mfg$-module homomorphism
$$
\pi_*:H^{r+*}_{\mfg,{cv}}(E,\E)\to H^*_\mfg(M,\F)
$$
in equivariant basic cohomology which we call {\em equivariant integration along the fibers}. In fact, if we regard an equivariant basic form $\omega\in \Omega^{*}_{cv}(E,\E)$ as an equivariant polynomial map $\omega:\mfg\to \Omega^{*}_{cv}(E,\E)$, then $\pi_*\omega$ is the equivariant polynomal map $\mfg \to \Omega^*(M,\F)$ given by $(\pi_*\omega)(X)=(\pi_*\omega)(X)$.

Now let $(M,\F)$ be a Killing foliation with structural Killing algebra $\mfa$ (see Section \ref{sec:Molino}). In this case
$$
\Omega_\mfa(M,\F)=S(\mfa^*)\otimes \Omega(M,\F)
$$
because $S(\mfa^*)^\mfa=S(\mfa^*)$, since $\mfa$ is abelian, and $\Omega(M,\F)^\mfa=\Omega(M,\F)$ by \cite[Lemma 3.15]{GT2010}. We will now consider the situation of Corollary \ref{cor:basicThomstratum}, where $N$ is a closed stratum, e.g.~a component of $C$, of codimension $r$. The natural foliation $\whF$ on the normal bundle $\nu N$ is also a Killing foliation with structural Killing algebra $\mfa$. With respect to these transverse actions the foliated bundle $\pi:(\nu N,\whF)\to (N,\F)$ is $\mfa$-equivariant, compare with Section \ref{sec:transverseaction}, inducing an $S(\mfa^*)$-module homomorphism, the {\em equivariant integration along the fibers}
$$
\pi_*:H^{r+*}_{\mfa,{cv}}(\nu N,\whF)\to H^*_\mfa(N,\F).
$$
This map is an isomorphism, since $\pi_*:H^{r+*}_{cv}(\nu N,\whF)\to H^*(N,\F)$ is an isomorphism by Corollary \ref{cor:basicThomstratum}, and because of \cite[Theorem 3.4]{McC} (or \cite[Theorem 6.7.1]{GS}) applied to the spectral sequences of $\Omega_\mfa(N,F)$ and $\Omega_{\mfa,cv}(\nu N,\whF)$ of \cite[Theorem 3.23]{GT2010}. Therefore its inverse $\tau_\mfa:H^*(N,\F)\to H^{r+*}_{cv}(\nu N,\whF)$, the {\em equivariant basic Thom isomorphism}, is given as in the ordinary case (see Section \ref{sec:BasicThom}) by
$$
\tau_\mfa(\omega)=\pi^*\omega\wedge \Phi_\mfa,
$$
where $\Phi_\mfa\in H_{\mfa,cv}^r(\nu N,\whF)$ is the unique class with $\pi_*\Phi_\mfa=1$. The equivariant basic Thom homomorphism can be realized on the level of forms with the use of the universal Thom form found by Mathai and Quillen (\cite[Equation (40)]{Meinrenken}). It is explicitly given in \cite[(40)]{Meinrenken} which has to be slightly modified to be of compact vertical support as in remarked in the same paragraph. The universal Thom form of the $\mathfrak{so}(r)$-action on $\RR^r$ then gives an equivariant basic Thom form $\Phi_\mfa$ by adapting \cite[Section 10.2]{GS} to the foliated principal bundle $(P,\E)$ of transverse oriented frames of $(\nu N,\whF)$ with its $\mfa$-invariant Riemannian basic connection form.

The composition of $\tau_\mfa$ with the $S(\mfa^*)$-homo\-mor\-phism $H^{*}_{\mfa,cv}(\nu N,\whF)\to H^{*}_\mfa(N,\F)$ induced by an inclusion of $(\nu N,\whF)\hookrightarrow (M,\F)$, where $\nu N$ is identified with a saturated tubular neighborhood of $N$ in $M$, will be denoted by
$$
i_*^N:H_\mfa(N,\F)\to H_{\mfa}(M,\F)
$$
We will call this map the {\em equivariant basic Thom homomorphism} of $N$ in $M$ which by construction extends the basic Thom homomorphism of $N$. As it will turn out later, this map is injective, which is not necessarily true for its nonequivariant counterpart.

\section{Localization of equivariant basic characteristic classes}\label{sec:localization}
\subsection{Basic characteristic classes}\label{sec:BasicCharacteristicClasses}
Let $\pi:(P,\E)\to (M,\F)$ be a foliated $\mathbf O(r)$- respectively $\mathbf{SO}(r)$-bundle with a basic connection form. Then the curvature form $\Omega$ is basic with values in ${\rm Ad} P$. Let $S(\mathfrak{so}(r)^*)^{\mathbf O(r)}$ respectively $S(\mathfrak{so}(r)^*)^{\mathbf SO(r)}$ be the algebra of $\mathbf O(r)$-invariant respectively $\mathbf SO(r)$-invariant polynomials on $\soq$. The latter has one additional generator, the Pfaffian $e$. For $p\in S(\mathfrak{so}(r)^*)^{\mathbf O(r)}$ the Pontryagin forms  $p(\Omega)$ and, in the case of $\mathbf{SO}(r)$ in addition the Euler form $e(\Omega)$, are basic. The corresponding basic cohomology classes in $H^*(M,\F)$ are independent of the choice of the basic connection and will be called {\em basic characteristic classes}. We define the basic Chern-Weil homomorphism
$$
{\rm cw}_\F:S^\bullet(\mathfrak{so}(r)^*)^{\mathbf O(r)}\to H^{2\bullet}(M,\F) \quad \mbox{resp.} \quad {\rm cw}_\F:S^\bullet(\mathfrak{so}(r)^*)^{\mathbf{SO}(r)}\to H^{2\bullet}(M,\F)
$$
by ${\rm cw}_\F(p)=[p(\Omega)]$.
The image is called the {\em basic (Euler-)Pontryagin ring} associated to $(P,\E)$. For a taut, transversely oriented Riemannian foliation $\F$ of codimension $q$ on a complete manifold $M$ such that $M/\oF$ is compact we can now define characteristic numbers. We take the Molino bundle $(\whM,\whF)$ with the canonical Riemannian basic connection as our foliated $\SOq$-bundle and choose a $p(\Omega)$ from the associated basic Euler-Pontryagin ring of degree $q$. Then $\int_\F p(\Omega)$ is called the corresponding {\em basic characteristic number}. Since the integral depends on the choice of a volume element $\nu_\mfa$ of the strucural Killing algebra, the individual numbers strictly speaking are not invariants of the Riemannian foliations. Nevertheless their ratios or the set of these numbers up to scaling are invariants.

\subsection{Localization}
Let $\F$ be a transversely oriented Killing foliation of codimension $q$ with strucural Killing algebra $\mfa$ on a complete manifold $M$ such that $M/\oF$ is compact. Under these assumptions we will localize top basic Pontryagin classes to the union $C$ of all closed leaves, and every top basic form in the equivariantly formal case.

We will see $\RR$ for the moment as an $\mfa$-dga, where the derivations $d, \iota_X$ and $\cL_X,\ X\in\mfa$ are trivial. Then $\int_\F:\Omega(M,\F)\to\RR$ is an $\mfa$-dga-homomorphism inducing a map $\int_\F:\Omega_\mfa(M,\F)\to S(\mfa^*)$, which is in fact given by $(\int_\F\omega)(X)=\int_\F\omega(X)$, for which the diagram
$$
\xymatrix{
  H^*(M,\F) \ar[d]^{\int_\F}& H_\mfa^*(M,\F)\ar[l]_{{\rm ev}_0} \ar[d]^{\int_\F} \\
  \RR & \ar[l]_{{\rm ev}_0} S(\mfa^*) }
$$
is commutative first on the level of forms; here the horizontal maps mean inserting zero into the equivariant form regarded as a polynomial map $\mfa\to \Omega(M,\F)$ respectively $\mfa\to\RR$. Passing to (equivariant) basic cohomology we obtain the above diagram. Note that the top map is surjective if and only if the action is equivariantly formal (\cite[Corollary 3.28]{GT2010}). This means that in general not every closed basic form has an extension as a closed equivariant basic form. On the other hand a basic characteristic class of the normal bundle of $(M,\F)$ can be extended in analogy to \cite{BerlineVergne} as we will see now. Let $p\in S(\soq^*)^\SOq$be an invariant polynomial. Furthermore let $\omega$ be the connection form of the foliated $\SOq$-bundle $\rho:P\to M$, the Molino bundle, with respect to the canonical Riemannian basic connection. Since the transverse action of $\mfa$ is isometric, $\omega$ is $\mfa$-invariant. Let $\Omega$ be the curvature form of $\omega$. We obtain a basic characteristic class $[p(\Omega)]\in H^*(M,\F)$. For each $X\in \mfa$ we define $L_X:P\to\soq$ by $L_X=-\iota_{X^*}\omega$. Then the polynomial map $p(\Omega+L_X)$ in $X$ is closed with respect to $d_\mfa$, i.e.~it is a closed equivariant basic form, compare \cite[Proposition 2.13(iii)]{BerlineVergne}. Its class in $H^*_\mfa(M,\F)$ is called an {\em equivariant basic characteristic class}. Clearly under the top map in the above diagram it is mapped to $p(\Omega)$.

We will now derive an integration formula for the foliation setting by adapting the treatment of the classical case in \cite[p.~8f]{AtiyahBott}. Let $i:(C,\F) \to (M,\F)$ be the inclusion of the union of all closed leaves in $M$. Let $\{C_i\}$ be the set of components of $C$. For the equivariant basic Thom map $i_*:H_\mfa(C,\F)\to H_\mfa(M,\F)$ we have $i^*i_*1=e_\mfa(\nu C,\whF)$,  where the latter is the equivariant basic Euler class. In other words the equivariant basic Thom class restricts to the equivariant basic Euler class. This follows from the shape of the universal Thom form and \cite[Theorem 8.2]{Meinrenken}. Since the equivariant basic Euler class $e_\mfa(\nu C,\whF)$ is a not a zero-divisor in $H^*_\mfa(C,\F)=S(\mfa^*)\otimes H^*(C,\F)$, $i_*$ is injective and induces an isomorphism in the localized module $\widehat H_\mfa(M,\F)=Q(\mfa^*)\otimes_{S(\mfa^*)} H_\mfa(M,\F)$; here $Q(\mfa^*)$ is the field of fractions of $S(\mfa^*)$. Then
$$
S=\sum_i\frac{i^*_{C_i}}{e_\mfa(\nu C_i,\whF)}
$$
is inverse to $i_*$ in $\widehat H_\mfa(M,\F)$. Thus for any equivariant basic class we have
$$
\omega=i_*S\omega=\sum_i\frac{i_*^{C_i}i^*_{C_i}\omega}{e_\mfa(\nu C_i,\whF)}.
$$
after localization. Applying $\int_\F:H_\mfa(M,\F)\to S(\mfa^*)$ to both sides of this equation  and using Lemma \ref{lem:restrictionintegral}, we obtain the next theorem with the constants $c_i$ from Section \ref{sec:PoincareDuality}.
\begin{thm}[ABBV-type localization formula for Killing foliations]
Let $\F$ be a transversely oriented Killing foliation on a complete Riemannian manifold $M$ such that $M/\oF$ is compact.
For any $\omega\in H_\mfa(M,\F)$ we have
$$
\int_\F\omega=\sum_i c_i \int_{C_i/\F}\left(\frac{i^*_{C_i}\omega}{e_\mfa(\nu C_i,\whF)}\right).
$$
\end{thm}
Note that we have a polynomial in $\mfa$ on the left, whereas the right side is a sum of rational functions. This formula allows us to compute basic characteristic numbers. Let $p(\Omega)$ be a basic characteristic class in $H^q(M,\F)$. Then $p(\Omega+L_X)\in H^*_\mfa(M,\F)$, the  corresponding equivariant basic characteristic class, is its equivariantly closed extension. We take this for $\omega$ in the localization formula. Inserting zero in both sides gives a formula for the basic characteristic number $\int_\F p(\Omega)$.
This computation is particularly simple if the closed leaves are isolated. We write $C=\bigcup L_i$. Then the right side can be expressed in terms of the weights of the isotropy representation of $\mfa$ on $\nu_{x_i} L_i$ for arbitrary choices of $x_i\in L_i$. Let $\Delta_i$ be the correponding set of weights.
\begin{cor}
Let $\F$ be a transversely oriented Killing foliation of even codimension $q$ on a complete Riemannian manifold $M$ such that $M/\oF$ is compact. Moreover, assume that the closed leaves $\{L_i\}$ are isolated. Then
$$
\int_\F p(\Omega)=\sum_i (-1)^{q/2}c_i \frac{p(L_X)}{\prod_{\alpha\in\Delta_i} \alpha(X)}.
$$
\end{cor}

\appendix
\section{Existence of a good saturated cover}\label{sec:cover}
In Definition \ref{def:finitetype} we defined what a good saturated cover is, and when a foliated manifold is of finite type. In this section we want to prove the following result.
\begin{prop}
 A Riemannian foliation $\F$ on a manifold $M$ that is complete with respect to an adapted bundle-like metric has a good saturated cover. If $M/\oF$ is compact, $(M,\F)$ is of finite type. Moreover, for any good saturated cover there is a basic partition of unity subordinate to it.
\end{prop}
See Definition \ref{def:finitetype} for the definition of a good saturated cover and foliations of finite type.
Via the correspondence principle described at the of Section \ref{sec:subMolino} the first statement amounts to finding a good invariant cover (see below) of the $K$-manifold $W$, where $W$ is the basic manifold ocurring in Molino's structure theory, see Section \ref{sec:subMolino} and $K$ is either $\SOq$ or $\Oq$, depending on whether $\F$ is transversely orientable or not. For the second statement in the proposition it then suffices to find a subordinate $K$-invariant partition of unity; this is stated in \cite[Theorem 1.45]{AlBiPe}.
In \cite[Section 1]{AlNi} a proof of the existence of a particular good invariant cover is given. Unfortunately that section contains a few mistakes. For a corrected proof it is helpful to discuss this briefly. First \cite[Remark 1.2.(ii)]{AlNi} is wrong for an orbit $Ky$ of the same dimension as $Kx$ but of smaller orbit type; a well-known $S^1$-action on the Klein bottle provides a counterexample. Consequently the map mentioned in the first sentence of the last paragraph of the proof of \cite[Theorem 1.3]{AlNi} is not necessarily surjective.  It is not possible to fix this by choosing $Ky$ to be of minimal orbit type in $U_1\cap U_2$ instead of just of minimal dimension. The problem is that $U_1\cap U_2$ might contain two different locally minimal orbit types. A good example is the action of the torus $T^2\subset \mathbf{SO}(5)$ on $S^4\subset\RR^5$ and $U_1$ and $U_2$ are balls around the two fixed points with nonzero intersection (if we flatten the metric on the intersection $U_1\cap U_2$ the $U_i$ are regular neighborhoods). In that case the $T$ in the quoted text does not necessarily intersect both orbit types transversally (as in the example of the torus action) and therefore the aforementioned map is not necessarily a surjective diffeomorphism.

Let $K$ be a Lie group acting properly on a complete manifold $W$ endowed with a $K$-invariant Riemannian metric. We show the existence of a good invariant cover.
\begin{dfn} A cover $\cU=\{U_\alpha\}$ of a $K$-manifold $M$ by open invariant sets $U_\alpha$ is called a {\em good invariant cover} if each nonempty finite intersection $V=U_{\alpha_0}\cap \cdots \cap U_{\alpha_p}$ is $K$-{\em contractible}, meaning that there is an orbit $Kx\subset V$ and a homotopy $H:V\times [0,1]\to V$ such that $kH(x,t)=H(kx,t)$ for all $k\in K, x\in V, t\in [0,1]$ (in other words $H_t$ is $K$-equivariant), and $H_0={\rm id}_V$ and $H_1(V)= Kx_{\alpha_i}$.
We call $(M,K)$ of {\em finite type} if it admits a finite good invariant cover.
\end{dfn}

A geodesic is called {\em transnormal} if it is orthogonal to every orbit it meets. If a geodesic is orthogonal to an orbit in one point, it is transnormal. We say a $K$-invariant set $U\subset M$ is {\em strongly transversely convex} if for any two orbits in $\overline U$ there is a unique (up to $G$-translation) minimal geodesic between them, whose interior lies completely in $U$. Note that this geodesic is orthogonal to both orbits by the variation formula and consequently transnormal.
For a $K$-orbit $N$ and $r\geq 0$ let $U_r(N)$ be the open tubular neighborhood of $N$ with radius $r$.
Then $N$ has positive injectivity radius $i(N)$, i.e., $\exp^\perp_N:\nu^rN\to U_r(N)$ is a diffeomorphism for any $0<r<i(N)$. Let $\rho:\overline U_r(N)\to N$ be the orthogonal projection.
\begin{lemma}\label{lem:tangency}
Let $N$ be a $K$-orbit. Then there is a real $0<c<i(N)$ such that for any $0<r<c$ and any transnormal geodesic $\gamma$ with $q:=\gamma(0)\in\partial U_r(N)$ and $\dot\gamma(0)\in T_q\partial U_r(N)$ remains outside $U_r(N)$ for some time.

Moreover, each tubular neighborhood of $Kx$ of radius $r$ with $0<r<c$ is strongly transversely convex.
\end{lemma}
\begin{proof}
The proof of the first part is similar to that of \cite[Lemma 3.4.1]{doCarmo}. We will only comment on the differences. The function $F(t,v)=|(\exp^\perp_N)^{-1}(\gamma_v(t))|^2$ measures the square of the distance of a geodesic $\gamma_v$ with initial unit vector $v\in TU_r(N), r<i(N)$ to $N$ instead of to $p$. In order to show $\frac{\partial F}{\partial t}(0,v)=0$ for $v\in T\partial U_r(N)$ we need to apply the Gauss Lemma for submanifolds (e.g.~\cite[Ch.~II, (4.14)(ii)]{Sakai}) instead of the ordinary Gauss Lemma. Now as in the reference one sees $\frac{\partial^2F}{\partial t^2}(0,v)=2|v|^2=2$ for a unit normal vector of $N$. By continuity there is a $0<c<i(N)$ such that $(0,v)$ is a minimum of $F$ for any $0<r<c$ and any unit vector $v\in T\partial U_r(N)$ belonging to a transnormal geodesic. This finishes the first part.

The proof of the second part is similar to that of \cite[Proposition 3.4.2]{doCarmo} relying on the first statement of this lemma. We consider a minimal geodesic $\gamma$ between two orbits in $U_r(Kx)$. It is transnormal and the tangent vector of $\gamma$ at maximal distance to $Kx$ fulfills the assumptions of the first statement of the lemma. One concludes the lemma as in \cite[Lemma 3.4.1]{doCarmo}.
\end{proof}
For each orbit $Kx$ let $r_x>0$ be such that $2r_x$ is smaller than the constant $c$ from Lemma \ref{lem:tangency}. In particular
\begin{enumerate}
	\item $\exp^\perp_{Kx}:\nu^{2r_x}Kx\to U_{2r_x}(Kx)$ is a diffeomorphism,
	\item $U_{r}(Kx)$ is strongly transversely convex for all $r$ with $0<r<2r_x$.
\end{enumerate}
Clearly this yields an invariant cover of $K$-contractible neighborhoods $U_x:=U_{r_x}(Kx)$. Let $\cU=\{U_{x_\alpha}\}$ be this cover or a subcover of it. Our aim is to show that any nonempty finite intersection $U_{x_{\alpha_1}}\cap \ldots \cap U_{x_{\alpha_k}}$ is also $K$-contractible. Fix one such intersection $U:=U_{x_{\alpha_1}}\cap \ldots \cap U_{x_{\alpha_k}}\neq\emptyset$ and let $x_i:=x_{\alpha_i}$. Let $\cW$ be the stratification of $W$ by $K$-components of orbit type manifolds (\cite[Section 3]{HT1}), and the stratum of $\cW$ containing a point $x\in W$ is denoted by $\cW_x$. We can assume $r_{x_1}\geq \cdots \geq r_{x_k}$. Then $x_{i+1}\in U_{2r_{x_i}}(Kx_i)$ by the triangle inequality. Property (1) implies that for any $y\in U_{2r_{x_i}}(Kx_i)$ the stratum $\cW_y$ contains the minimal geodesic segment from $Kx_{i}$ to $y$, possibly minus the starting point. It follows $x_i\in \overline\cW_{x_{i+1}}$ and therefore $\cW_{x_i}\subset\overline\cW_{x_{i+1}}$.

Next we want to show that $\cW_{x_{i}}\cap U_{x_{1}}\cap \ldots \cap U_{x_{i}}\neq\emptyset$ for all $i$. Fix $i$. Choose $y\in U_{x_{1}}\cap \ldots \cap U_{x_{{i}}}$. Let $\gamma$ be a minimal geodesic from $\overline\cW_{x_{i}}$ to $y$, whose starting point we denote by $y'$. We have $d(y',Kx_j)\leq d(y',y)+d(y,Kx_j)\leq 2r_{x_j}$ for all $j\leq i$, since $d(y',y)=d(y,\overline\cW_{x_i})\leq d(y,Kx_i)\leq r_{x_i}\leq r_{x_j}$. So $\gamma$ lies in $U_{2r_{x_j}}(Kx_j)$ by property (2) and starts orthogonally to the submanifold $\overline\cW_{x_i}\cap U_{2r_{x_i}}(Kx_i)=\cW_{x_i}\cap U_{2r_{x_i}}(Kx_i)$ by the variation formula. Thus $\gamma$ is transnormal. We want to see that $d(y',Kx_j)\leq d(y,Kx_j)<r_{x_j}$ for all $1\leq j\leq i$ which implies $y'\in \cW_{x_i}\cap U_{x_{1}}\cap \ldots \cap U_{x_i}$ and we are done. Fix $j$. The transnormal geodesic $\gamma$ starts tangentially to $\partial U_{d(y',Kx_j)}(Kx_j)$ and is outside $U_{d(y',Kx_j)}(Kx_j)$ at least for small $t$ by the first part of Lemma \ref{lem:tangency} and our choice of $r_{x_j}$. Assume $d(y',Kx_j)> d(y,Kx_j)$ then $y\in U_{d(y',Kx_j)}$. Since $U_{d(y',Kx_j)}$ is strongly equivariantly convex it contains the interior of $\gamma$, contradiction.

Now let $i_0$ be the smallest number $i\in\{1,\ldots,k\}$ for which $\cW_{x_i}\cap U_{x_{1}}\cap \ldots \cap U_{x_{k}}\neq\emptyset$. Then $\cW_{x_{i_0}}$ is the only locally minimal orbit type in $U$ (see (\cite[Definition 3.4]{HT1})).
Since the intersection of strongly transversely convex invariant sets is itself strongly transversely convex, and $U$ is $K$-contractible as we will see in Lemma \ref{lem:intersectionset}, $\cU$ will be a good saturated cover, concluding the proof of the Proposition.

Before that we need the following lemma.
\begin{lemma}\label{lem:isotropygeodesic}
Let $\gamma:[a,b]\to W$ be a minimal geodesic between two $K$-orbits. Then the isotropy group $K_{\gamma(t)}$ is the same for all $t\in (a,b)$. Moreover, we have $K_{\gamma(a)}\supset K_{\gamma(t)}\subset K_{\gamma(b)}$ for all $t\in (a,b)$.
\end{lemma}
\begin{proof}
Clearly, the geodesic $\gamma$ is transnormal. Let $t\in [a,b]$.  By the linearity of the isotropy representation on $\nu_{\gamma(t)}K\gamma(t)$ the isotropy groups $K_{\gamma(s)}$ are the same for all $s\neq t$ in a neighborhood of $t$ in $[a,b]$. Thus the isotropy groups $K_{\gamma(t)}, t\in [a,b]$ are the same except possibly in a finite set. Assume there is such an exception at  $t\in (a,b)$. Let $t$ be the first such point. Then $K_{\gamma(t)}\supset K_{\gamma(s)}$ for all $s\neq t$ in a neighborhood. Take $k\in K_{\gamma(t)}\backslash K_{\gamma(s)}$. Then $k\gamma|[a,t)$ is another segment from $K\gamma(a)$ to $\gamma(t)$ different from $\gamma|[a,t)$. The concatenation of this segment with $\gamma|[t,b]$ has the same length as $\gamma$. One can shorten this curve at the angle at $\gamma(t)$. This contradicts the minimality of $\gamma$. Thus $K_{\gamma(t)}$ is the same for all $t\in (a,b)$.
\end{proof}

Now let $U$ be an arbitrary strongly transversely convex invariant neighborhood of an orbit $Kx$. Let $V=V(U,Kx)\subset\nu Kx$ be the set of $\dot\gamma(0)$ such that $\gamma$ is a minimal geodesic parametrized by arc length between $Kx$ and an orbit in $U$. Clearly $V$ is $K$-invariant and $\exp^\perp_{Kx}:V\to U$ is a surjective $K$-equivariant map that is not a diffeomorphism in general.
\begin{lemma}\label{lem:intersectionset}
Let $U$ be an open strongly transversely convex invariant set containing only one locally minimal orbit type and let $Kx$ be an orbit of that orbit type in $U$. Then $\eta:=\exp^\perp_{Kx}:V(U,Kx)\to U$ is a diffeomorphism and in particular $U$ is $K$-contractible.
\end{lemma}
\begin{proof}
We have already remarked that $\eta$ is surjective, and we want to show injectivity. Assume $\eta(v_1)=\eta(v_2)=:y$ for some $v_1,v_2\in V=V(U,Kx)$ with foot points $x_1$ respectively $x_2$ in $Kx$. Since $U$ is strongly transversely convex there is a $k\in K$ such that $v_1=k_*v_2$ so $k\in K_y$. From Lemma \ref{lem:isotropygeodesic} we know that the isotropy groups along the geodesic $\gamma:[0,\|v_2\|]\to U$ with $\|v_2\|\dot\gamma(0)=v_2$ are constant except possibly at the end points $x_2$ and $y$. We denote the stratification by $K$-components of orbit type manifolds in $U$ by $\cW$. Let $Ky'$ be of locally minimal orbit type in the closure $\overline\cW_y$. By assumption there is only one locally minimal orbit type in $U$, namely $\cW_x$, so we have $\cW_x=\cW_{y'}\subset\overline\cW_y$. Thus $\gamma|(0,\|v_2\|]$ lies in $\cW_y$ and $K_y=K_{\gamma(t)}$ for all $t\in (0,\|v_2\|]$. So $k\in K_{\gamma(t)}$ for all $t\in [0,\|v_2\|]$ by Lemma \ref{lem:isotropygeodesic} which implies $v_1=v_2$. This proves that $\eta$ is injective.

We will now show that $\eta$ has full rank. Assume that $d\eta$ is singular at $v\in V$, i.e., $\eta(v)$ is a focal point of $Kx$ and let $\sigma:[0,\infty)\to M$ be the geodesic with $\dot\sigma(0)=v/\|v\|$. But then $\sigma|[0,t]$ is not a minimal geodesic between $Kx$ and $K\sigma(t)$ for $t>\|v\|$. This contradicts the injectivity of $\eta$. Therefore $\eta$ is a diffeomorphism.
\end{proof}


\end{document}